\documentclass[10pt,a4paper,times]{article}
\usepackage[latin1,applemac]{inputenc}
\usepackage{amsmath, amssymb, amsfonts, amsthm, amscd, epsfig, multicol, marvosym, float,mathrsfs,enumerate,appendix,float,color,upgreek}
\usepackage[left=2cm,top=3cm,bottom=3cm,right=2cm]{geometry}
\usepackage{multirow} 
\usepackage[mathscr]{euscript} 

\numberwithin{equation}{section}

\begin{document}
\title{Motions about a fixed point by hypergeometric functions: new non-complex analytical solutions and integration of the herpolhode}
\author{Giovanni Mingari Scarpello \footnote{giovannimingari@yahoo.it}
\and Daniele Ritelli \footnote{Dipartimento scienze statistiche, via
Belle Arti, 41 40126 Bologna Italy, daniele.ritelli@unibo.it}}
\date{}
\maketitle
\begin{abstract}
The present study highlights  the dynamics of a body moving about a fixed point, and provides  analytical closed form solutions. Firstly, for the symmetrical heavy body, that is the Lagrange-Poisson case, we compute the second (precession, $\psi$) and third (spin, $\varphi$) Euler angles in explicit and real form by means of multiple hypergeometric (Lauricella) functions. Secondly, releasing the weight assumption but adding the complication of the asymmetry,  by means of elliptic integrals of third kind, we provide the precession angle $\psi$ completing the treatment of the Euler-Poinsot case. Thirdly, by integrating the relevant differential equation, we reach the finite polar equation of a special motion trajectory named the {\it herpolhode}. Finally, we keep the symmetry of the first problem, but without weight, and take into account a viscous dissipation. The use of motion first integrals -adopted for the first two problems- is no longer practicable in this situation, therefore the Euler equations, faced directly, are driving to particular occurrences of Bessel functions of order $-1/2$.

\noindent\textcolor{blue}{This is a pre-print of an article published in Celestial Mechanics and Dynamical Astronomy. The final authenticated version is available online at: DOI: 10.1007/s10569-018-9837-5}
\vspace{0.5cm}


\noindent {\sc msc classification:} 	70E17, 33C65, 33C10, 33E05
\end{abstract}
\section{Introduction}
\subsection{Aim of the paper}
The motion of a rigid body about a fixed point has been a subject of mathematical studies and generalizations for almost three centuries. In this article we want to achieve four goals.

\begin{enumerate}

\item For the symmetrical heavy solid (case Lagrange-Poisson), we calculate the second and third Euler angles in explicit form by means of multiple hypergeometric functions, thus obtaining, via the classic rotation matrix, the position of a gyroscope's arbitrary point as seen by the fixed axes. 

\item For the asymmetrical body free of torques (case Euler-Poinsot), we compute the precession during time by means of elliptic integrals of third kind. Integrating then the general form of the relevant nonlinear differential equation, we get the finite polar equation of a special trajectory named the {\it herpolhode}.

\item The third problem we analyze, providing exact solutions, concerns a body with symmetry as in the first problem, but without the weight and adding dissipation to it. In fact air resistance acts as an obvious cause of energy loss for the motion of the top. In the majority of the real situations, exact solutions are not available: this is due to the lack of first integrals.

\item  We select eighty items among papers, books and handbooks and  provide an up-to-date state of the art describing some topics-as the Serret-Andoyer variables-not so well known besides the astronomers.  
\end{enumerate}

 Modern treatments \cite{law} sometimes use the elliptic functions of Weierstrass (WEF) like $\wp(z; g_2, g_3)$ or the non-elliptic ones like $\sigma$ or $\zeta$: but they are hardly effective to computation. In fact physicists are less familiar with the WEF, say with $\wp(z; g_2, g_3)$, partly because the form of its fundamental period parallelogram depends on the invariants $g_2$ and $g_3$ and on the modular discriminant $\Delta = g_2^{3} - 27 g_3^{2}$ as well. 
In this way, the applications of WEF are usually analyzed in terms of four different cases based on the signs of $(g_3,\Delta)$, namely $ [(-,-),(-,+),(+,-),(+,+)]$, with two special cases $(g_3=0,\Delta=0)$ and $(g_3=0,\Delta>0)$. This is in contrast with the fundamental period parallelogram of the Jacobi elliptic functions which remains rectangular for all values of $m \ne 1$ where $m$ is a function of the roots $(e_1, e_2, e_3)$ of the cubic at right hand side of the ODE: ${\rm d}s/{\rm d}z=4s^3-g_2 s-g_3$ defining $\wp$. The greatest advantage of the WEF, however, could be that the function $\wp (z)$ itself and its derivative $\wp'$ can be used to represent any doubly-periodic function as $f (z) = A(\wp) \wp' + B(\wp)$, where $A$ and $B$ are arbitrary functions.

As example of a different approach, in \cite{Pal} the motion about a fixed point is studied of a body under the effect of a force whose potential depends on an angle $\theta$ alone. Furthermore the body is under the effect of another potential changing with $\theta$ as well: but the prolonged use of imaginary unit and of unusual elliptic and nearly elliptic functions, weakens the work, and makes useless this and many other papers which could be of some use even today. 
Some antique treatments as Jacobi \cite{Jac1} up to Greenhill \cite{gree}, Appell \cite{appe2} and Whittaker \cite{whit} have set their calculations on {\it complex variables} which are difficult to be implemented: anyway this point will be further analyzed throughout this paper.

We stayed in the {\it real field}: we will build explicit solutions through Jacobi elliptic functions cn, sn, dn, elliptic integrals of first and third kind, Appell function and those Lauricella hypergeometric functions that have been profitably used by cosmologists \cite{Kr} and by ourselves \cite{msr}. These multiple\footnote{We use, according to Exton \cite{exx} the term \lq\lq multiple'' meaning that the hypergeometric  is a function of two, or more, variables.} hypergeometric functions are exposed in several special treatises as \cite{exx}, \cite{slat} and \cite{bailey}.

\subsection{Nomenclature}

The main symbols recurring through this paper are:

\begin{tabular}{ l l }
  $A,\, B,\, C $ & body principal moments of inertia  \\
  $D$ & fictitious body moment of inertia  \\
  $\mu$ &  constant of viscous torques  \\
  $m$ &  homogenizing constant   \\
  $\delta,\,a,\,b,\,c,\lambda $ &  parameters about viscous torques  \\
  $\boldsymbol{\omega}$ & instantaneous angular speed vector\\
  $p,\, q,\, r$& projections of  instantaneous vector $\boldsymbol{\omega}$ on the body Cartesian coordinate system\\
  $p_M,\, q_M$& maximum values of $p(t),\, q(t)$: problem of torqueless top\\
  $P,\, Q,\, R$ & projections of  instantaneous vector $\boldsymbol{\omega}$ on the space Cartesian coordinate system\\
  $x,\, y,\, z$ & coordinates of a body point with respect to the body Cartesian coordinate system\\

  $X,\, Y,\, Z$ & coordinates of a body point with respect to the space Cartesian coordinate system\\
  $\theta,\,\varphi,\,\psi$ & Euler angles: nutation, intrinsic rotation=''spin'', precession\\
  $x_G,\,y_G,\,z_G$ & coordinates of the centre $G$ of gravity with respect to the body  Cartesian coordinate system\\

  $X_G,\,Y_G,\,Z_G$ & coordinates of the centre $G$ of gravity with respect to the space  Cartesian coordinate system\\
  $\alpha_1,\,\alpha_2,\,\alpha_3$ & direction cosines of the body axes respect to the $OX$ axis of space Cartesian coordinate system\\
 $\beta_1,\,\beta_2,\,\beta_3$ & direction cosines of the body axes respect to the $OY$ axis of space Cartesian coordinate system\\

 $\gamma_1,\,\gamma_2,\,\gamma_3$ & direction cosines of the body axes respect to the $OZ$ axis of space Cartesian coordinate system\\
${\bf i},\, {\bf j},\, {\bf k}$ & unity vectors of the body Cartesian coordinate system\\
${\bf I},\,{\bf J},\,{\bf K}$ & unity vectors of the space Cartesian coordinate system\\
${\bf v_G}$ & centre $G$ of gravity velocity \\
$O',\, O$ & origins of the body and space Cartesian coordinate systems\\
$M\boldsymbol{g}$ & weight of the rigid body\\
$E_0$ & total constant energy\\
${\bf K}_O$ & moment of momentum taken with respect to $O$\\
$K_{z,0}$ & initial value of the $z$-component of ${\bf K}_O$\\
${\bf K^0}_O$ & initial value of ${\bf K}_O$\\
$s$ & traditional shortening for $\cos\theta(t)$\\
$s(0),\, s_1,\, s_2,\, s_3$ & specific determinations of  $\cos\theta(t)$\\

\end{tabular}

\begin{tabular}{ l l }

${\bf K}(k)$ & complete elliptic integral of first kind of modulus $k$\\
$F(\varphi,k)$ & incomplete elliptic integral of first kind of modulus $k$ and amplitude $\varphi$\\
$\Pi(\varphi, \alpha, k) $ & incomplete elliptic integral of third kind of amplitude $\varphi$, parameter $\alpha$, modulus $k$\\
$J_\nu$ & $\nu$-order Bessel functions of the first kind\\
$Y_\nu$ & $\nu$-order Bessel functions of the second kind\\
$t,\, \tau\, $ & time, adimensional time\\
${{\rm F}_{D}}^{(n)}$ & $n$-variable multiple hypergeometric Lauricella function\\
 $\mathcal {G}$ & constant appearing in the herpolhode equation\\
$\rho_1,\,\rho_2,\,\rho $ & the herpolhode's vector radius determinations\\
$ \chi(\rho) $ & the herpolhode's polar anomaly corresponding to $\rho$\\
$\varepsilon,\, \sigma $ &  parameters of third kind elliptic integrals appearing in  precession, torqueless case\\
{$a,\, b,\, c $} & functions of $A,\, B,\, C,\, D $  homogeneous to a reciprocal of a moment of inertia\\
\end{tabular}

\vspace{0.3cm}

Some other few symbols are used throughout the text as ease variables without a specific meaning and do not need to be listed here.

\subsection{Historical minimal outlook}

The motion of a rigid body about a fixed point became part of the mathematical literature thanks to J. d'Alembert (1717--1783) in 1749, see \cite{Alem1}. The following year, L. Euler (1707--1783) in \cite{Eul1} tried to look at such motions with a more general point of view, but with the handicap of moments of inertia changing during the motion. Both mathematicians were lacking of the principal axes of inertia and their properties. The first rigorous treatment of these axes was due to J. A. Segner (1704--1777) who proved, in 1755, \cite{Seg}, that free rotation is possible through a minimum of three individual axes, being more than three for special cases of symmetry (spheres, etc.). Euler acknowledged the strength of Segner's reasoning, and was the first to see that these axes all had to pass through the centre of gravity.
As a consequence, d'Alembert, 1761, \cite{Alem2}, took advantage by simplifying his formul\ae\      and Euler issued his {\it  Theoria motus corporum solidorum seu rigidorum} \cite{Eul2}, 1765, in which the angles $\theta$, $\phi$ and $\psi$ are introduced. But what Euler had obtained after lengthy calculations, J. L. Lagrange (1736-1813) by his own method, assuming $\theta,\,\varphi,\,\psi$ as free coordinates, accomplished ({\it M\'{e}canique Analytique}, 1788) in a few lines \cite{Lagra}. S.D. Poisson (1781-1840) gave in his {\it Trait\'{e} de M\'{e}canique} \cite{Pois}, 1833,  an elementary proof of the Euler equations without introducing the Euler angles providing also the system \eqref{Poi} of the {\it kinetic equations} named after him. A.M. Legendre (1752-1833) is the first author of the idea of substituting, as a means of investigation, an ideal ellipsoid, having certain relations with the actually revolving body. Although he conducts his own investigations on principles that are different, he seems well aware of the use which might be made of this happy conception, as one can read in his {\it Trait\'{e}} \cite{Leg}, tome I page 410.

The problem of motion of a $M\bf{g}$-heavy rigid body about a fixed point under the sole gravity can be reduced \cite{Eul2} to find the solutions of the following system of three ordinary differential equations in the unknown $p,\, q,\, r$, functions of time, namely the projections of the instantaneous vector $\boldsymbol{\omega}$ on the body Cartesian coordinate system:
\begin{equation}\label{Eu}
\begin{cases}
A{\dot{p}}=(B-C)qr+Mg\left(y_{G}\gamma_{3}-z_{G}\gamma_{2}\right)\\
B{\dot{q}}=(C-A)rp+Mg\left(z_{G}\gamma_{1}-x_{G}\gamma_{3}\right)\\
C{\dot{r}}=(A-B)pq+Mg\left(x_{G}\gamma_{2}-y_{G}\gamma_{1}\right)
\end{cases}
\end{equation}
where  $A,\, B,\, C $ are the body moments of inertia with respect to its principal axes and  $x_{G},\,y_{G},\, z_{G}$ the (constant) coordinates of the centre of mass seen by the body Cartesian coordinate system. Finally, the $\gamma$'s are unknown functions of time providing the direction cosines, as seen by the body Cartesian coordinate system, of the vertical axis along which the gravity is acting. They are linked to $p,\, q,\, r$ by the Poisson \cite{Pois} kinetic equations:
\begin{equation}\label{Poi}
\begin{cases}
{\dot{\gamma_{1}}}=r\gamma_{2}-q\gamma_{3}\\
{\dot{\gamma_{2}}}=p\gamma_{3}-r\gamma_{1}\\
{\dot{\gamma_{3}}}=q\gamma_{1}-p\gamma_{2}
\end{cases}
\end{equation}
Perhaps inspired by the mentioned Legendre's idea, A. L. Cauchy (1789-1857) \cite{Cau} introduced, in 1827, his inertia or {\it momental ellipsoid}. Anyway, in 1851, the seminal {\it Th\'eorie nouvelle de la rotation des corps} by L.
Poinsot (1777-1859) \cite{Poins}  appeared with a geometric interpretation of the rigid body motion about a fixed point, which allows to represent it as clearly as that of a moving point. The ellipsoid rolls without slipping on an invariable plane orthogonal to the moment of momentum vector; while  $\boldsymbol{\omega}$'s endpoint signs on it a curve polhode which rolls on  the curve herpolhode signed by itself on the fixed plane.
 We will do some effort in integrating the herpolhode's differential equation. Coming again to the Euler-Poisson differential system \eqref{Eu} and \eqref{Poi} despite to its simple appearance, only a few exact solutions have been obtained until now. It consists of six nonlinear equations of the first order and will depend on six constants but, due to the relationship among the cosines, there will be just five of them. The system \eqref{Eu} will be solvable whenever admits sufficient number of  {\it first integrals}: the reduction of order is then applied because each first integral eliminates one equation. For the reduction of the problem to quadratures it is then sufficient to have only four independent first integrals which do not contain time. Three classic algebraic integrals are known for the heavy body: the energy integral, the conservation law of the angular momentum along the vertical direction, and the geometrical constraint on the direction cosines, namely:
\begin{equation}\label{trio}
\begin{cases}
Ap^2+Bq^2+Cr^2+2Mg\left(x_{G}\gamma_{1}+y_{G}\gamma_{2}+z_{G}\gamma_{3}\right)=E_{0}\\
Ap\gamma_{1}+ Bq\gamma_{2}+Cr\gamma_{3}=K_{0}\\
\gamma_{1}^2+\gamma_{2}^2+\gamma_{3}^2=1
\end{cases}
\end{equation}
For a body free of torques, the first integral of angular momentum will imply the constancy of the whole norm of  momentum. A fourth algebraic first integral for arbitrary values of the coefficients of the six equations is then necessary, but it is unknown a priori. The literature on the subject is immense: the book of Leimanis, 1965, \cite{Leim} quotes about six hundred entries. So far, the following basic cases of integrability have been found for the heavy body.

\noindent $\bullet$ {\bf Case Euler-Poinsot}: 

The assumptions are $x_{G}=y_{G}= z_{G}=0$ so that the system \eqref{Eu} becomes homogeneous and the fourth first integral is:
\[
Ap^2+Bq^2+Cr^2=\text{constant}.
\]
It is the case of {\bf body free of torques} {also} studied by L. Poinsot. In absence of  loading couples, for instance, when the body is rotating about its center of gravity, the equation can be driven to quadratures even without {in lack of} special body symmetry. C.G.J. Jacobi (1804-1851) in \cite{Jac1} and \cite{Jac2}, both of 1849, applied the elliptic functions in order to evaluate the Euler angles and the components $p,\, q,\, r$ of the angular speed.

\noindent $\bullet$ {\bf Case Lagrange-Poisson}:

The assumptions are $x_{G}=y_{G}=0;\, z_{G}\ne 0$ and  $A=B \neq C$, so that the fourth first integral is: 
\[
r(t)=r_{0}=\text{constant}.
\]
The body is under the effect of its weight and the constraint is assumed to be smooth. Throughout this paper we {are also committed in} facing some problems with the integration of such a heavy rotor.

\noindent $\bullet$ {\bf Case of complete kinetic symmetry}:

The  assumptions are $A=B=C$ so that the fourth first integral is: 
\[
x_{G} p+y_{G} q+z_{G}r=\text{constant.}
\]

\noindent $\bullet$  {\bf Case Koval\`{e}vskaja}:

Such a case was solved \cite{Kov} by S. V.  Koval\`{e}vskaja (1850-1891) in 1888 by assuming
that two of the principal moments of inertia are equal and double the third, $A=B=2C$, and that its $G$ is situated in the plane of the equal moments of inertia, $z_{G}=0$.
So that the fourth first integral she found is:
\[
\left[C\left(p^2-q^2\right)+Mgx_{G}\gamma_{1}\right]^2+\left[2Cpq-Mgx_{G}\gamma_{2}\right]^2=\text{constant}.
\]
Her memoir stimulated a large number of investigations regarding new particular solutions of the general problem of the new case. She also proved that in the general case, the fundamental system of six equations does admit single-valued solutions containing five arbitrary constants and having no singular points other than poles only in the aforementioned four cases. Furthermore, A. M. Lyapunov (1857-1918) proved that only for the above four basic cases, the six functions: $p,\, q,\, r$ and the $\gamma's$ are single-valued ones for any arbitrary initial values.

\noindent $\bullet$ {\bf Post-Lyapunov developments}:

Beyond the four basic cases of integrability, if further restrictions are imposed on the constant values of the energy and of the angular momentum about the vertical, and on $A,\, B,\, C,$ and on the $G$ position with respect to co-mobile axes, then a solution becomes possible. Having a system of $n$ ordinary differential equations of the first order in normal form, let its solution be ${\bf x}:$ we say that the real function $E$ is a first integral of a continuous dynamic system over $W\subset\mathbb{R}^n$ if $E$ is a differentiable function over $W,$ such that:
\[
E({\bf x}(t))=E({\bf x} (0))=E({\bf x}_0)
\]
for each solution, namely for each initial condition ${\bf x}(0)$. If such a generality fails, we have what Levi-Civita calls \cite{levi} {\it partial first integrals} namely true only for special starting velocity field \lq\lq atto di moto''. As far as we are concerned, so far there are {ten} of such restricted cases of integrability:

\begin{enumerate}
\item { W. Hess \cite{Hess} considered in 1890 the case in which: }
\begin{enumerate}
\item {the mass centre $G$ of the rigid solid is on the normal at the fixed point O driven to one of the planes of circular section of the ellipsoid of inertia, assumed not of rotation;} 
\item {at the initial time, the moment of momentum is situated in the respective plane of circular section.} 
\end{enumerate}
{This case has been furthermore studied again by G.G. Appelrot, P.A. Nekrasov, B.K. Mlodzeevski, N.E. Jukovski, S.A. Chaplygin, R. Liouville and others, due to its importance. T. Manacorda \cite{Man} considered in 1950 the case in which the moment of the
given external forces with respect to the fixed point is normal to the straight line $OG,$
concluding that the moment of momentum  must have the same property and then confirming the considerations of Hess and Appelrot.}

{If, in particular, at the beginning, the moment of momentum is horizontal,
then the constant of areas vanishes and $G$ moves as a 
mathematical pendulum; the plane of the circular section is rotating about a fixed horizontal straight line, normal to the trajectory of the centre $G.$ The trajectory of the point $P(0,\sqrt{B/M},0),$ end point of the mean axis of the ellipsoid of inertia on the sphere centered  in $O$ of radius $\sqrt{B/M}$ will be a loxodrome. If the mass centre has an oscillatory motion, then the point $P$ oscillates
along an arc of the loxodrome and if the centre $G$ has an asymptotic motion, then the motion of the rigid solid tends asymptotically to a rotation about the aforementioned axis of
inertia. The rigid solid is, in this case, named {\it loxodromic pendulum}.
}

\item {The Bobylev-Steklov Case. }

{In this case, which has been independently  put in evidence by D. Bobylev and V.A. Steklov in 1896, one assumes that $A = 2B,$ the centre of mass being on the $Ox_1$-axis. In such a case the Euler differential system of angular speed components becomes tractable and the problem is completely solved in terms of elliptic integrals.}

\item {In 1894, O. Staude \cite{Stau} and B.K. Mlodzeevski have shown independently that the equations of Euler and Poisson allow a simple infinity of solutions if no one restriction is imposed about the moments of inertia and about the position of the mass center of the body. These solutions correspond to rotations of the body about axes fixed in the body and in space.
Mlodzeevski showed that such {\it permanent axes of rotation} can be the principal axes of inertia, if these axes are horizontal (the case of the physical pendulum) or a family of vertical axes; in the latter case, studied in detail by Staude, the magnitude of the instantaneous angular velocity is constant.
}

\item {Second case of Steklov. If $B>A>2C$ with the fixed point at the origin. Then permanent rotations about axes lying in the plane $z=0$ are stable, while rotations about those in the plane $y=0$ are unstable if the square of cosine of the angle between the permanent axis of rotation and the $x$-axis is greater than $1/3$}

\item {Goriatchev-Chaplygin.
Goriatchev showed that when $A=B=4C$ the mass center of the body is in the plane of equal moments of inertia at $O$ and the projection of the body angular momentum on the vertical through $O$ is zero, then a fourth particular algebraic integral exists and the solution depends on three particular constants. Chaplygin improved next year such a result and obtained a solution containing four arbitrary constants.}

\item {Second Goriatchev.
At the beginning of XX century, Goryachev, Steklov and Chaplygin put the problem to find new cases of integrability, where the first integrals, independent of time, be algebraic, not containing arbitrary constants, see e.g. Teodorescu \cite{teodorescu}. Conditions under which such integrals exist impose restrictions upon the body principal moments of inertia either by equality or inequality restricting the values of the ratios of two principal moments of inertia.}

\item {Second Chaplygin, 1904. Making special assumptions, see Leimanis, \cite{Leim}, on ratios $C/A$ and $B/A,$ he obtains time as a function of two elliptic integrals of the variable $u$ generated by the Hermite transformation
\[
 \frac{1+b^2c^2z^4}{z^2}=u 
\]
where $z^3=p(t).$ Inverting the first kind elliptic integral, $u$ can be found as a function of time; the angular speed $p$ will come by solving with respect to $p$ the equation
\[
 b^2 c^2 p^{2/3}+p^{-2/3}=u.
\] 
}

\item {N. Kovalevski \cite{Ko} considered in 1908 the case in which the mass centre is on a principal axis, searching all the cases in which the components of angular speed can be expressed in the form of polynomials of the third degree in $p.$ Thus, he found again the particular cases studied by Goryachev, Steklov and Chaplygin, as well as a new case where the special relationship
$AC=9(B-C)(A-2B)$ between the principal inertia moments $A,\, B,\, C$ holds.
}

\item {Grioli \cite{Grioli} considered, 1947, a heavy rigid body free to rotate about a fixed point $O$ different from its mass center $G$ and with the central ellipsoid of inertia asymmetrical, that is a non revolution ellipsoid. He showed that among the possible motions of such a heavy asymmetric body there are $\infty^2$ regular non-degenerate precessions having as the figure axis through to $G$ straight line orthogonal to the plane of one of the circular cross sections of the ellipsoid itself. His proof was based on purely dynamical considerations without integrating the equations what was done later by Gulyaev (1955).}

\item {Cases of spherically symmetric body and axis-symmetric body subjected to time-constant body-fixed torque, by Romano, 2008 \cite{romano}. The Romano's paper analyzes several cases. By his formulae 8, 9, 10, the author by $p_0, q_0$ and U constructs the complex variables z and $\nu$ (by (10)). The Kummer function $_1F_1$ has then to be evaluated at complex values of parameters and of the variable. At the end he arrives at $G(z,c)$ which will provide the stereographic complex rotation. 
Further notices on very recent developments will be given in the following paragraph.
}

\end{enumerate}

Details on the relevant contents can be found in \cite{Bogo} and in \cite{Leim} with a bibliography on the Russian original papers. R. Liouville \cite{Liou}, 1896, tried to find all cases in which the system of six equations admits a fourth  first integral algebraic and independent of time. At the end, P. Burgatti (1868-1938) proved, 1910,  \cite{Burg} that each fourth algebraic integral not dependent on time always turns out to be a combination of other three, except for the four fundamental cases of Euler, Lagrange, spherical and Koval\`{e}vskaja.

Recently, the motion of a rigid body about a fixed point under the influence of a Newtonian force field has been investigated \cite{NES} obtaining the necessary and sufficient conditions for some functions to be a fourth first integral of the governing equations.
We have to comment two papers.

Vrbik \cite{vrbik} uses the Lagrange technique to find the corresponding ODE for the three Euler angles. He first analyzes a steady solution ignoring nutation successively added as a periodic perturbation of a small amplitude. Next he passes to an \lq\lq exact'' treatment but no solution is obtained, providing only some lines of programming in Mathematica$^{\text{\textregistered}}$ for the precession but not for nutation and spin.

\section{The $\boldsymbol{F}$-loaded symmetrical body}

Let us consider a homogeneous rigid body whose centre of mass is $G$ and referred to a space Cartesian coordinate system $OXYZ,\,OZ$ upwards, and to a body Cartesian coordinate system $O'xyz$, with $O$ and $O'$ overlapped. If we assume all the loads being exerted on it to reduce to a force ${\bf F}=(F_x, F_y, F_z)$ applied in $G,$ then the general expression of the axial torque generalizing the right hand side of the first of \eqref{Eu} becomes:
\[
F_x( y_G \alpha_3-\alpha_2 z _G)+F_y( z_G \beta_1-\beta_3 x _G)+ F_z( y_G \gamma_3-\gamma_2 z _G)
\]
and so on for the other ones. The unknowns of the problem are then twelve: three body-components $p,\,q,\,r$ of $\boldsymbol{\omega}$  plus the nine variables $\alpha_j,\, \beta_j,\, \gamma_j,\, j=1,\,2,\,3 $, direction cosines of the body axes $Oxyz$ with respect to the space Cartesian coordinate system. The equations are twelve: the three Euler joint to nine kinetic Poisson equation (namely \eqref{Poi} plus three analogous for $\beta$'s  and those for $\alpha$'s) linking the cosines to the $p,\,q,\,r$ components.
The problem is intractable at such a level of generality, so that we will consider a particular case of the highest interest.

\section{The heavy symmetrical body (Lagrange-Poisson)}

The problem of the rigid heavy body rotating about a fixed point assumes that in $G$ only the weight vector is applied. The body is then moving under the simultaneous effects of the weight and of the constraint reaction which, due to the smoothness of the support plane, is vertical too. As a consequence, the moments cardinal equation projected on the body Cartesian coordinate system axes provides the system \eqref{Eu}. Due to our assumptions on $G,\, x_G=y_G=0$ and $z_G\neq0$, we are then facing the problem under a double non-zero (and variable) component of torque, as we can read on \eqref{Eu}. 
We get:
\begin{equation}\label{Eu-weight}
\begin{cases}
A{\dot{p}}=(B-C)qr+Mg\left(-z_{G}\gamma_{2}\right)\\
B{\dot{q}}=(C-A)rp+Mg\left(z_{G}\gamma_{1}\right)\\
C{\dot{r}}=(A-B)pq
\end{cases}
\end{equation}

To \eqref{Eu-weight}, one has to append the Poisson kinematic equations \eqref{Poi}: so that the system is reduced to only six ordinary differential equations. 

All the problem data are: four body\footnote{We assume that the spinning top mass $M$ does not change during its motion.} specifications $A,\ C,\ M,\ z_G$ and six initial values: $\varphi_0,\, \theta_0,\, \psi_0,\, r_0,$ $K_{z,0},\  E_0$.

Let us start with the existence of three first integrals. The lack of any dissipation provides us the energy as the first integral:
\begin{equation}\label{Ene}
\frac{1}{2}(Ap^2+Bq^2+Cr^2)+MgZ_G=E_0 
\end{equation}
being the vertical $OZ$ fixed axis oriented upwards.
Furthermore, the external (weight {plus} reaction) couple's moment being perpendicular to $OZ,$ its  projection on $OZ$ will be zero, so that the  component $K_{z,0}$ of the $O$ moment of momentum ${\bf K}_O$ will be constant.
Projecting it on the space Cartesian coordinate system, one will provide the second first integral:
\begin{equation}\label{Momm}
 Ap\gamma_{1}+ Bq\gamma_{2}+Cr\gamma_{3}=K_{z,0}.  
\end{equation}
A third first integral will stem from the geometrical link of cosines of the body axes with respect of the $OZ$ space axis:
\begin{equation}\label{Coss}
 \gamma_{1}^2+\gamma_{2}^2+\gamma_{3}^2=1 .
\end{equation}
The knowledge of three first integrals \eqref{Ene}, \eqref{Momm} and \eqref{Coss} is not enough to solve the six ordinary differential equations system by quadratures. It is well known, see for instance \cite{toma}, page 401, that even if we have six equations,  it is sufficient to know four independent prime integrals, {namely those generating a full rank Jacobian matrix}, first integrals of the considered differential system. In fact, having the last Jacobi multiplier value 1, a fifth integral can be determined as independent of the other ones, so that the system can be solved by knowing four first integrals. Accordingly, the treatment of the heavy body requires specific assumptions in order to obtain such a fourth first integral. Under certain restrictions concerning the $G$-location and the values of the moments of inertia $A,\, B,\, C $ towards co-mobile axes $O'x,\,O'y,\,O'z$, then such a fourth integral can be found for any initial values of energy and momentum.

Each solid body whose central inertia ellipsoid {is of revolution}, is named a gyroscope, such that the $O'z$ gyroscopic axis holds the centre of gravity: in such a case the ellipsoid will be {of revolution} with respect to any other point on the axis.
{Therefore we can infer that} $A=B,$ while the latter implies that $x_G =y_G =0,\, z_G>0$. In such a way the third of \eqref{Eu} provides:
\begin{equation}\label{Quarto}
r(t)=r_{0}={\rm constant}.
\end{equation}
Therefore, as a consequence of the assumption on the mass distribution we see that the angular speed around the $O'z$ axis remains unchanged in time. In such a way the first two integrals become:
\begin{equation}\label {Bibus}
\begin{cases}
\dfrac{1}{2}A(p^2+q^2)+Mgz_{G}\gamma_{3}=E_{0}-\dfrac{1}{2}Cr_0^2\\
\\
A(p\gamma_{1}+ q\gamma_{2})+Cr_0\gamma_{3}=K_{z,0}\\
\end{cases}
\end{equation}
There are three possible analytical paths for the problem of the body around a fixed point. The first, less used due to analytic difficulties, is to work directly on the differential system \eqref{Eu}. The second does not use the Euler equations at all but the Lagrange ones, assuming the Euler angles as {attitude parameters}. For instance Pars's treatise \cite{Pars}  shows this, but after the nutation resolvent, he does not go on with the analysis on spin and precession but deviates in analyzing steady motions.

Finally, one can work on the first integrals and we will follow this way for the heavy body and torque-free body problems.
\subsection{The nutation angle $\boldsymbol{\theta}$}
We refer briefly to the classic inquiry about the nutation $\theta$ as a basis for all our  further developments. Starting from \eqref{Bibus}, by working aside on the expressions on the left hand sides, we obtain $p^ 2 + q ^2= \hat {\rho}^2 (-\cos\theta +h-c \lambda^2)$ and afterwards:
\begin{equation}\label {zaia}
 p\gamma_{1}+ q\gamma_{2}=  \hat {\rho} (\hat{k}-c \lambda \cos\theta) 
 \end{equation}
which will give:
\[
\left[(\cos\theta)^.\right]^2=(p^ 2 + q ^2)(1-\cos^2\theta)-(p\gamma_{1}+ q\gamma_{2})^2
\]
so that, taking into account \eqref{Coss} and putting: 
\[
s=\cos\theta (t); \quad  \hat {\rho}^2 =2\frac{Mgz_{G}}{A}; \quad \ h=\frac{E_0}{Mgz_{G}}; \quad c=\frac{C}{A};  \quad\ \hat {k}=\frac{K_{z,0}}{A \hat {\rho}};\quad \ \lambda =\frac{r_0}{ \hat {\rho}},
\]
at the end we have:
\[
\dot{s}^2=f(s)=\hat {\rho}^2 \left[(h-c\lambda^2)-s\right](1-s^2)-\hat {\rho}^2 \left[\hat{k}-c\lambda s\right]^2.
\]
The above expression  is commonly referred to as the {\it resolvent of the heavy gyroscope}.

In the literature it is possible to see how many possible cases have been studied, but we are interested only in that of {\it major physical sense}, namely founded on:
\begin{multicols}{2}
\begin{description}
\item[ \bf {Assumption 1}]  $c\lambda \ne \mp \hat{k}$;
\item[ \bf {Assumption 2}]  $\ f(\mp1)=-(\mp c \lambda -\hat{k}^2)<0.$
\item[ \bf {Assumption 3}]  $f (s(0)) > 0.$
\end{description}
\end{multicols}
 Accordingly, let $s(0)$ be the value of $\cos\theta (t)$ for $t=0.$ Applying to $f(s)$ the roots existence theorem between $s=-1$ and  $s=s(0)$ a simple root will exist between them, say  $s_1<0$. Again, between $s=s(0)$ and $s=1$ there exists a further simple root, say  $s_2>0$. The cubic $f(s)$ has then two real roots lying between $-1$ and $+1$; whence its third real root, say $s_3>1$, does not have the intrinsic meaning of the previous ones, we can then provide the resolvent as the elliptic curve:
 \begin{equation} \label {essepunto}
  \dot{s}^2=f(s)= \hat{\rho}^2(s-s_1)(s-s_2)(s-s_3)
    \end{equation} 
which has real meaning only if $s_1<s<s_2.$ We can therefore see how the nutation angle $\theta$ changes periodically during time: starting with the value $\arccos (s(0))$, it oscillates, always bounded between $\theta_1=\arccos (s_1)$ and  $\theta _2=\arccos (s_2)$.
 Time will then be provided through a incomplete elliptic integral of first kind. 
 For these and further integrations we will make reference to the Appendix. By the  integral  $I_7$ in the Appendix, we infer:
\[
\int_{s(0)}^{\cos\theta}\frac{{\rm d}s}{\sqrt{f(s)}}=\frac{2}{\sqrt{s_2-s_1}}\left[F\left(\arcsin\sqrt{\frac{\cos\theta-s_1}{s_2-s_1}},k\right)-F\left(\arcsin\sqrt{\frac{s(0)-s_1}{s_2-s_1}},k\right)\right]
\]
being
\[
k=\sqrt{\frac{s_2-s_1}{s_3-s_1}} 
\]
the modulus of the Jacobi sine-amplitude function. Carrying the inversion:
\[
L=\hat{\rho}t+\frac{2}{\sqrt{s_3-s_1}}\,F\left(\arcsin\sqrt{\frac{s(0)-s_1}{s_2-s_1}},k\right),
\]
so that nutation succeeds to be:
\begin{equation} 
\label{sn}
s(t)=\cos\theta(t)=s_1+\left(s_2-s_1\right){\rm sn}^2\left\{\frac{\sqrt{s_3-s_1}}{2} \left[\hat{\rho}t+\frac{2}{\sqrt{s_3-s_1}}\,F\left(\arcsin\sqrt{\frac{s(0)-s_1}{s_2-s_1}},k\right)\right],k\right\}
\end{equation}
such a formula is well known, \cite {land}, \cite{synge}, \cite{wittenburg}. Lurie, \cite{lurie} at page 380 of his {\it Analytical Mechanics} treatise, writes the ODE of the nutation, but for its further developments refers to a monograph of  V.V. Golubev \cite{Golubev} which does not contain our results explicitly. A detailed discussion of the resolvent equation is founded on the inquiry of the real roots of the cubic. We did three assumptions leading to three real simple roots so that, if $s(0)$ is, {for example}, between $s_1$ and $s_2$, a $\theta$-periodic motion takes place. It attracted our interest for its major physical sense: as a matter of fact the body's figure axis will move in the space describing a conic surface limited by two round cones of vertical axis and half-aperture $\theta_1=\arccos(s_1)$ and  $\theta_2=\arccos(s_2)$. This is called the {\it nutation movement of the gyroscope's axis}.

\subsection{The precession angle $\boldsymbol{\psi}$}

After that the eulerian angle of nutation $\theta$ has been obtained explicitly as a periodic function\footnote{The periods of elliptic function $\rm{sn}$ are $4\boldsymbol{K}(k)$ and $2i\boldsymbol{K}'$. The periods of $\rm{cn}$ are $4\boldsymbol{K}$ and $2\boldsymbol{K}+2i\boldsymbol{K}'$ being $\boldsymbol{K}(k)$ the complete elliptic integral of first kind and $k'$ the complementary modulus. The imaginary period of such functions is meaningless to the motion. Nevertheless Appell \cite{appe1} shows that in the Poinsot torque-free motion, under convenient initial conditions one can associate the body movements twice so that the real period of the former equates the latter imaginary period divided by $i=\sqrt{-1}$.
} of time, let us pass to the other angles $\varphi $ and $\psi $. We spent some room to nutation (treated also e.g. in \cite{synge}) because we are building on it our original formulae for ${\varphi}$ and ${\psi}$. 
On this purpose, let us recall the components $p,\,q,\,r$ of the instantaneous angular speed can be expressed in terms of the Euler angles and their time-derivatives thanks to a set of formul\ae:
\begin{equation}\label{pqr}
\begin{cases}
p=\dot{\theta}\cos\varphi+ \dot{\psi} \sin\theta \sin\varphi\\
q=-\dot{\theta}\sin\varphi+ \dot{\psi} \sin\theta \cos\varphi\\
r=\dot{\varphi}+ \dot{\psi} \cos\theta
\end{cases}
\end{equation}
By inserting \eqref{pqr} and \eqref{Poi} in the second of \eqref{Bibus} we obtain: $K_{z,0}-Cr_0 \cos\theta=\dot{\psi} \sin^2\theta$ so that:
\begin{equation}\label{psipunto}
\frac{\rm d \psi}{\rm d t}=\hat{\rho}\frac{k-c\lambda s}{1-s^2}.
\end {equation}
Having been determined $s=\cos\theta$ as elliptic function of time, the above integration could be tried in the time domain through the elliptic functions. We prefer a different approach: by eliminating time between \eqref{psipunto} and \eqref {essepunto}, we arrive at:
\begin{equation}
\label{iper1}
\frac{\rm d \psi}{\rm d s}=\frac{(1-\hat{c}\lambda s/\hat{k})}{(1-s^2)\sqrt{(s-s_1)(s-s_2)(s-s_3)}}
\end {equation}
 to be integrated with respect to $s$ when $\psi$ goes from $\psi_0$ to its generic value. Such an integration has been carried out in Appendix  as $I_1$ by means of Lauricella multiple hypergeometric functions ${\rm F}_D^{(4)}$  of four independent variables, each of them being a different modulation of  $\cos\theta$ and then of time. In such a way we get:
\begin{equation}
\label {psilaur}
\psi (\theta)-\psi_0=\frac{2\sqrt{\cos\theta-s_1}}{\left(1-s_1^2\right)
   \sqrt{(s_3-s_1) (s_2-s_1)}}\left[\left(1-\frac{\hat{c}\lambda}{\hat{k}}\,s_1\right)\,\Xi-\frac{\hat{c}\lambda}{3\hat{k}}\left(\cos\theta-s_1\right)\,\Lambda\right]
\end{equation}
where 
\begin{align}\label{chi}
\Xi&=\mathrm{F}_{D}^{(4)}\left( \left. 
\begin{array}{c}
\frac12;1,1,\frac12,\frac12 \\[2mm]
\frac32
\end{array}
\right| \frac{\cos\theta-s_1}{1-s_1},-\frac{\cos\theta-s_1}{1+s_1},-\frac{\cos\theta-s_1}{s_1-s_3},-\frac{\cos\theta-s_1}{s_1-s_2}\right)\\
\\
\Lambda&=\mathrm{F}_{D}^{(4)}\left( \left. 
\begin{array}{c}
\frac32;1,1,\frac12,\frac12 \\[2mm]
\frac52
\end{array}
\right| \frac{\cos\theta-s_1}{1-s_1},-\frac{\cos\theta-s_1}{1+s_1},-\frac{\cos\theta-s_1}{s_1-s_3},-\frac{\cos\theta-s_1}{s_1-s_2}\right).\label{lambda}
\end{align}
Nevertheless - for completeness sake - the integral \eqref{iper1} has been also evaluated in terms of elliptic functions, see the Appendix.

\subsection{The intrinsic rotation $\boldsymbol{\varphi}$ (spin)}
The last angle is the one around the $z$-gyroscopic axis, namely the intrinsic rotation or spin. We can start from the third of  \eqref{pqr}, which, by means of the fourth first integral \eqref{Quarto} becomes $\dot{\varphi}=r_0- \dot{\psi} \cos\theta.$ Having already found the time laws of $\cos\theta $ and $\psi$, we could proceed to integrate the above differential equation, which is extremely difficult.
Accordingly, we will do something else, putting \eqref{psipunto} into the previous one, so that:
\[
\frac{\rm d \varphi}{\rm d t}=r_0 - \hat{\rho}\hat{k}\frac{(1-\hat{c}\lambda s/\hat{k})s}{1-s^2}.
\]
Again, by eliminating time between this equation and \eqref {essepunto}, we arrive at:
\begin{equation}
\label{iper2}
\frac{\rm d \varphi}{\rm d s}=\frac{(1-\hat{c}\lambda s/\hat{k}) s}{(1-s^2)\sqrt{(s-s_1)(s-s_2)(s-s_3)}}
\end {equation}
 to be integrated with respect to $s$ when $\varphi$ goes from $\varphi_0$ to its generic value. At the Appendix we see how such an integral $I_2$ can be developed by means of the Lauricella's function ${\rm F}_D^{(4)}$ representation theorem, or, alternatively, by elliptic functions.
We get:
\begin{equation}
\begin{split}
\label {philaur}
\varphi (\theta)-\varphi_0=&\frac{2\sqrt{\cos\theta-s_1}}{(1-s_1^2)
   \sqrt{(s_3-s_1) (s_2-s_1)}}\times\\
&\times\left[\left(1-\frac{\hat{c}\lambda}{\hat{k}}\,s_1\right)\,\Xi+\frac{1-2\hat{c}\lambda/\hat{k}\,s_1}{3}\left(\cos\theta-s_1\right)\,\Lambda-\frac{\hat{c}\lambda /\hat{k}}{5}\left(\cos\theta-s_1\right)^2\,\Omega\right]
   \end{split}
\end{equation}
where $\Xi$ and $\Lambda$ are defined in $\eqref{chi}$ and \eqref{lambda} respectively, while $\Omega$ is provided by \eqref{omega} below:
\begin{equation}\label{omega}
\Omega=\mathrm{F}_{D}^{(4)}\left( \left. 
\begin{array}{c}
\frac52;1,1,\frac12,\frac12 \\[2mm]
\frac72
\end{array}
\right| \frac{\cos\theta-s_1}{1-s_1},-\frac{\cos\theta-s_1}{1+s_1},-\frac{\cos\theta-s_1}{s_1-s_3},-\frac{\cos\theta-s_1}{s_1-s_2}\right)
\end{equation}
Finally, injecting the time function $\cos\theta$ \eqref{essepunto} in both \eqref{psilaur} and \eqref{philaur}, we get eventually the $\varphi $ and $\psi$ time laws of precession and spin.

\subsection{The overall body motion with respect to the space Cartesian coordinate system}
In order to describe the overall motion of the rigid body, first of all we consider its orientation given by the direction cosines of the body axes (respect to the $OZ$ axis of space Cartesian coordinate system)
$\gamma_1,\,\gamma_2,\,\gamma_3$. They are linked to Euler angles by:
\begin{equation}\label{gammas}
\begin{cases}
\gamma_{3}= \cos\theta\\
\gamma_{2}=\sin\theta \cos\varphi\\
\gamma_{1}=\sin\theta \sin\varphi
\end{cases}
\end{equation}
so that, having found $\theta$ and $\varphi$ as functions of time, we can evaluate each of $\gamma'$s, even if, of course, $ \gamma_{1}$ will be quicker thanks to  \eqref{Coss}. Putting the so  computed functions $\gamma_{1},\gamma_{2},\gamma_{3}$ in \eqref{zaia} and in the third of \eqref{Poi}, we will be allowed to detect $p(t),\, q(t)$ as:
\begin{equation}\label{pq}
\begin{cases}
p(t)=\dfrac{\hat {\rho}(1-c\lambda \gamma_3)\gamma_1 -\gamma_2 \dot{\gamma_3}}{1-\gamma_3^2}\\
\\
q(t)=\dfrac{\hat {\rho}(1-c\lambda \gamma_3)\gamma_2 +\gamma_1 \dot{\gamma_3}}{1-\gamma_3^2}
\end{cases}
\end{equation}
In such a way our knowledge of body motion is almost complete as far as it concerns the body Cartesian coordinate system. In the spinning top we have three main points: its fixed point $O\equiv O'$ assumed as double origin of two Cartesian coordinate systems; the centre $G$ of gravity; the {\it apex} $A\equiv(0,0,1)$, namely the point at unity distance from $O$ and belonging to the mobile axis $O'z$.

\subsubsection {The apex $\boldsymbol{A}$}

Let us go to motion of the apex. 
Starting from the knowledge of Euler angles, we will pass to formulae providing the space components of $\boldsymbol{\omega}$, namely:
\begin{equation}\label{PQR}
\begin{cases}
P=\dot{\theta}\cos\psi+ \dot{\varphi} \sin\theta \sin\psi\\
Q=\dot{\theta}\sin\psi- \dot{\varphi} \sin\theta \cos\psi\\
R=\dot{\psi}+ \dot{\varphi} \cos\theta
\end{cases}
\end{equation}
The apex space position can be found by means of the matrix rotation starting from its body Cartesian coordinate system coordinates. We have:
\begin{equation}\label{matt}
\left( \begin{matrix}
X_A \\
 Y_A\\
 Z_A 
\end{matrix}\right)=\left( \begin{matrix}
\alpha_1&\alpha_2&\alpha_3 \\
 \beta_1&\beta_2&\beta_3\\
 \gamma_1&\gamma_2&\gamma_3 
\end{matrix}\right)\left( \begin{matrix}
 x_A \\
 y_A\\
 z_A 
\end{matrix}\right)
\end{equation}
Being $x_A=y_A=0,\, z_A=1$ and minding the cosines of $O'z$ with respect to space Cartesian coordinate system are given by:
\begin{equation*}
\begin{cases}
\alpha_{3}= \sin\theta \sin\psi\\
\beta_{3}=-\sin\theta \cos\psi\\
\gamma_{3}=\cos\theta 
\end{cases}
\end{equation*}
we get
\begin{equation*}
\left( \begin{matrix}
X_A \\
 Y_A\\
 Z_A 
\end{matrix}\right)=\left( \begin{matrix}
 \sin\theta \sin\psi\\
-\sin\theta \cos\psi\\
\cos\theta
\end{matrix}\right)
\end{equation*}
The basic formula of rigid kinematics provides:
\begin{equation*}
{\bf v}_A(t)=\boldsymbol{\omega}\wedge (A-O)=\det\left( \begin{matrix}
\textbf{I}& \textbf{J} & \textbf{K} \\
P(t) & Q(t) & R(t)\\
X_A & Y_A& Z_A 
\end{matrix}\right)
\end{equation*}
So that, having previously computed all the six functions inside the matrix, the parametric equations for the apex are obtained and the problem is solved.

\subsubsection {Arbitrary point $\boldsymbol{Q}$}
We wish to know the motion of an arbitrary point, say $Q,$ whose body Cartesian coordinates  are given by:
\[
\left( \begin{matrix}
x_Q \\
y_Q\\
z_Q 
\end{matrix}\right)
\]
Such a (constant) vector has then to be put, instead of that concerning $A$, in formula \eqref{matt}. In order to perform the matrix product all nine cosines will be necessary, see for instance \cite{agos} p. 209, as functions of the Euler angles:
\begin{equation}\label{9coseni}
\begin{cases}
\alpha_1=\cos \psi \cos\varphi -\sin\psi  \\
\alpha_2=-cos\psi \sin\varphi- \sin\psi \cos\varphi \cos\theta \\
\alpha_3= \sin\theta \sin\psi\\
\beta_1=\sin\psi \cos\varphi +\cos\psi\sin\varphi\cos\theta\\
\beta_2=-\sin\psi\sin\varphi+ \cos\psi\cos\varphi\cos\theta\\
\beta_3=-\sin\theta\cos\psi\\
\gamma_1= \sin\theta \sin\varphi\\
\gamma_2=\sin\theta \cos\varphi\\
\gamma_3= \cos\theta
\end{cases}
\end{equation}
The matrix product will provide $X_Q(t),\,Y_Q(t),\,Z_Q(t)$. At the end, the usual rigid body formula will give ${\bf v}_Q(t).$ 
In such a way we provided a very compact and new hypergeometric treatment of the heavy rigid body  motion, a top whose symmetry (case Lagrange-Poisson) gave us the required fourth first integral.

\subsubsection {The centre $\boldsymbol{G}$ of gravity}
In order to follow the motion of the centre of gravity, we proceed as for the apex, of course putting $z_G$ instead of 1 in column vector of  right hand side of \eqref{matt}. In such a way we arrive at the vector
\[
\left( \begin{matrix}
X_G \\
Y_G\\
Z_G 
\end{matrix}\right)
\]
and as usually we compute  ${\bf v}_G$ whose determination has a peculiarity: by means of the first cardinal dynamics equation
\[
{\bf \Phi}+Mg {\bf K}=M\frac{{\rm d}{\bf v}_G}{{\rm d}t}
\]
we finally can compute the constraint reaction $ \bf \Phi$. 

\section{The asymmetrical body free of torques (Euler-Poinsot)}

The rigid body motion about a fixed point can also be solved without any special symmetry, i. e. $A>B>C$ but in such a case it must be free of torques and then either hung at its centre of gravity or in free fall. In such a case we talk of the Poinsot motion whose study had been started by Euler, as we know.

\subsection{Nutation and spin}
We recall the classic analysis adding our new treatment of the precession. Two first integrals of the Euler equations hold in the lack of dissipation, namely that of energy and the other one about the norm of moment of momentum. The body being fixed at its centre of gravity, we put there the origin of both superimposed Cartesian coordinate systems so that with  $x_{G}=y_{G}= z_{G}=0$ in \eqref{trio} we get:
\[
Ap^2+Bq^2+Cr^2=E_{0}.
\]
Furthermore, due to the lack of torques, the whole norm of ${\bf K}_O$ shall be constant in amount and its direction unchanged whenever the body moves. Being that its components on the principal axes are given by $Ap,\, Bq,\, Cr$,  we have:
\[
A^2 p^2+B^2q^2+C^2r^2=\|{\bf K}_O\|^2
\]
as one  could get from \eqref{Eu} by putting to zero the right hand sides, multiplying by $2Ap, 2Bq, 2Cr$, integrating and adding. 
On the contrary, the instantaneous direction of the rotational speed $\boldsymbol{\omega}$  is not stationary and its magnitude 
$\|\boldsymbol{\omega}\|=\sqrt{p^2+q^2+r^2}$
 also will change during time.
The equations \eqref{Eu} for an asymmetrical body free of torques become:
\begin{equation*}
\begin{cases}
\dfrac{{\rm d}p}{{\rm d}t}=\dfrac{B-C}{A}qr\\
\\
\dfrac{{\rm d}q}{{\rm d}t}=\dfrac{C-A}{B}pr\\
\\
\dfrac{{\rm d}r}{{\rm d}t}=\dfrac{A-B}{C}pq
\end{cases}
\end{equation*}
Comparing them with the system of the Jacobian elliptic functions:
\begin{equation*}
\begin{cases}
\dfrac{{\rm d}}{{\rm d}u}{\rm cn}u=-{\rm sn}u\  {\rm d}nu\\
\\
\dfrac{{\rm d}}{{\rm d}u}{\rm sn}u={\rm cn}u\ {\rm d}nu\\
\\
 \dfrac{{\rm d}}{{\rm d}u}{\rm dn}u=- k^2 {\rm sn}u\ {\rm cn}u
\end{cases}
\end{equation*}
 such a close likeness shall mean some role of these functions in describing the body free motion.
A modern treatment of the so-called Poinsot case proceeds on the tracks of Jacobi who, in 1849, provided explicit expressions to $p,\,q,\,r$ by means of elliptic functions of time. In summary, the problem {input} data are: body specifications $A,\ B, \ C,\ M,\ z_G$ and the initial values: $q_0,\, \psi_0,\, \|{\bf K}_O\|,\,  E_0$.

Landau \cite{land}, solves as:
\begin{equation}\label{pqerre}
\begin{cases}
p(t)=\sqrt{\dfrac{2E_0 C-\|{\bf K}_O\|^2}{A(C-A)}}{\rm cn}(\tau,\hat{k})=p_M{\rm cn}(\tau,\hat{k})\\
\\
q(t)=\sqrt{\dfrac{2E_0 C-\|{\bf K}_O\|^2}{B(C-B)}}{\rm sn}(\tau,\hat{k})=q_M {\rm sn}(\tau,\hat{k})\\
\\
r(t)=\sqrt{\dfrac{\|{\bf K}_O\|^2-2E_0 A}{C(C-A)}}{\rm dn}(\tau,\hat{k})
\end{cases}
\end{equation}
under the assumption $\|{\bf K}_O\|^2 >2E_0 B, q(0)=q_0=0$ and being 
\begin{equation}\label{tauandti}
\tau=t\sqrt{\dfrac{(C-B)(\|{\bf K}_O\|^2-2E_0 A)}{ABC}}
\end{equation}
a non-dimensional time,  and 
\begin{equation}\label{hatkappa}
\hat{k}=\sqrt{\dfrac{(B-C)(2E_0C-\|{\bf K}_O\|^2)}{(C-B)(\|{\bf K}_O\|^2-2E_0 A)}}   
\end{equation}
the elliptic modulus of the Jacobian elliptic functions. They are periodic and after one period of time, the vector $\boldsymbol{\omega}$  returns to its initial position relative to the body Cartesian coordinate system. Nevertheless, the top itself is not in the same position relative to the body Cartesian coordinate system.
In order to describe such a position at any time, it is necessary to compute {some attitude parameters, for example} the Euler angles.
For the purpose, we know the moment of momentum ${\bf K}_O$ does not change with respect to the space Cartesian coordinate system $OXYZ$, so that with $\gamma_i,\, i=1,\ 2,\ 3$ the direction cosines of the space axis $OZ$ with respect to the body Cartesian coordinate system $O'xyz$,we get:
\begin{equation}\label{ApBqCr}
\begin{cases}
Ap(t)=\|{\bf K}_O\|\gamma_1\\
\\
Bq(t)=\|{\bf K}_O\|\gamma_2\\
\\
Cr(t)=\|{\bf K}_O\|\gamma_3
\end{cases}
\end{equation}
and putting there the $\gamma$ functions provided by the last three from \eqref{9coseni} we obtain:
\begin{equation}
\begin{cases}
\cos \theta=\dfrac{Cr(t)}{\|{\bf K}_O\|}
\\[3mm]
\tan\varphi=\dfrac{A}{B} \dfrac{p(t)}{q(t)}
\end{cases}
\end{equation}
and then, finally:
\begin{equation}\label{2angoli}
\begin{cases}
\cos\theta(t)=\sqrt{\dfrac{C(\|{\bf K}_O\|^2-2E_0 A)}{\|{\bf K}_O\|^2 (C-A)}}{\rm dn}(\tau,\hat{k})\\
\\
\tan\varphi(t)=\sqrt{\dfrac{A(C-B)}{B(C-A)}}\dfrac{{\rm cn}(\tau,\hat{k})}{{\rm sn}(\tau,\hat{k})}\end{cases}
\end{equation}

\subsection {A short digression}

To the free top problem, namely the rotator with a fixed point but free of forces (Euler-Poinsot case) the solution of the angle $\varphi(t)$ is given (p. 151 of \cite{whit}) as:
\[
\exp(2i\varphi)=C\frac{\Theta_{01}(\frac{\lambda t -ia}{2K})}{\Theta_{01}(\frac{\lambda t +ia}{2K})}\exp\left\{\left[\frac{2id}{A}+\frac{\lambda}{K}\frac{\Theta_{01}' (ia/2K)}{\Theta_{01} (ia/2K)}\right]t\right\}
\]
in such a formula $C,\, d,\, k$ and $\lambda$ are constant values coming from the problem physical data and $a$ will stem from an algebraic relationship with data $A,\,B,\,c,\,d,\,D$. Symbols $\Theta_{01}$ and $\Theta_{01}'$ refer to a Theta series of complex argument and to its derivative.  
The complex argument requires some work. In order to compute K one has to evaluate the \lq\lq nome'' $q$ by means of the series
\[
q=\frac{1}{2}b^2+\frac{1}{16}b^{10}+\frac{15}{212}b^{18}+\cdots
\]
where $b=\tan(1-k^2)^{1/4}$. Furthermore:
\[
\sqrt{\frac{2K}{\pi}}=1+2q+2q^4+2q^9+\cdots
\]
Up to this point one can have the $\nu$ values
\[
\frac{\lambda t -ia}{2K}, \,\frac{\lambda t +ia}{2K}
\]
to put inside the Theta series
\[
\Theta_{01} (\nu)=1-2q\cos(2\pi\nu)+2q^4\cos(4\pi\nu)+\cdots
\]
and in its derivative with respect to $\nu$, namely $\Theta_{01}'={\rm d}\Theta_{01}/{\rm d}\nu$.

Comparing this long approach to formula \eqref{2angoli} bis, we feel confirmed in our criterion of sticking to real variables not following such a complex treatment, probably at that time not intended for a practical use.

\subsection{The  precession}

The third step will be to detect the third angle law $\psi (t)$. By eliminating $\dot{\theta}$ from the first two of \eqref{pqr} we get:
\[
\frac{\rm d \psi}{\rm d t}=\frac{p\sin\varphi+q\cos\varphi}{\sin\theta}.                                                                         \]
By inserting in the first two of \eqref{ApBqCr} their expressions \eqref{9coseni}, squaring and adding, we get the precession differential equation:
\[
\frac{\rm d \psi}{\rm d t}=\|{\bf K}_O\|^2\dfrac{A p^2+B q^2}{A^2 p^2+B^2q^2}
\]
where the right-hand side consists of all the elliptic functions above. Putting 
\[
\frac{B^2}{A^2}\frac{q_M^2}{p_M^2}=\delta^2,\ 1-\delta^2=\varepsilon^2,\ \gamma^2=\frac{\varepsilon^2}{1+\varepsilon^2}<1,
\]
we finally obtain:
\begin{equation}
\psi(\tau)=\psi_0+\dfrac{\|{\bf K}_O\|}{A(\delta^2+\varepsilon^2)}\int_{0}^{\tau}\frac{{\rm cn}^2(\tau,\hat{k})}{1-\gamma^2 {\rm sn}(\tau,\hat{k})}{\rm d}\tau+\dfrac{\|{\bf K}_O\|B\delta^2}{A^2}\int_{0}^{\tau}\frac{{\rm sn}^2(\tau,\hat{k})}{1-\varepsilon^2 {\rm sn}(\tau,\hat{k})}{\rm d}\tau,
\end{equation}
which requires the integration of rational functions of Jacobian elliptic functions. We have been  driven to elliptic integrals of third kind  as shown in the Appendix as $I_4$ and  and $I_3$. We obtain:
\begin{equation}
\label{psitorqueless}
\psi(\tau)=\psi_0+\dfrac{\|{\bf K}_O\|}{A(\delta^2+\varepsilon^2)}\mathcal {T}_1+\dfrac{\|{\bf K}_O\|B\delta^2}{A^2}\mathcal {T}_2
\end{equation}
where:
\begin{equation}\label{terzaspecie}
\begin{cases}
\mathcal {T}_1=\dfrac{1}{\varepsilon^2}\left[\Pi({\rm am}(\tau,\hat{k}),\varepsilon^2,\hat{k})-\tau\right]\\
\\
\mathcal {T}_2=\dfrac{1}{\gamma^2}\left[(\gamma^2-1)\Pi({\rm am}(\tau,\hat{k}),\gamma^2,\hat{k})+\tau\right]
\end{cases}
\end{equation}
Notice that an integral like that we solved in terms of $\Pi$ is met in \cite{synge} (first vol., form.61b, p. 257) but not solved.
The precession angle is then provided by \eqref{psitorqueless} as a function of non dimensional time $\tau$ through elliptic integrals of third kind  $\mathcal {T}_1$ and $\mathcal {T}_2$ of parameters $\varepsilon^2$ and $\gamma^2$ and whose argument is the amplitude of the elliptic integral of first kind. The common modulus of all of them is $\hat{k}$ given by \eqref{hatkappa}. At the end one will pass to the effective time through \eqref{tauandti}.
Of course the knowledge gained of all the Euler angles will allow to compute the cosines by means of \eqref{gammas}. The classic Euler-Poinsot case admits a class of simple but non-trivial integrable generalizations which modify the Poisson equations describing the motion of the body in space, see \cite{fedorov2013}.

\subsection{Modern activity about the Euler-Poinsot  problem}

Although the Euler-Poinsot formulation via the Euler angles is quite known among engineers, nevertheless astronomers prefer the canonical variables as capable of a better modeling e.g. of the disturbing torques, of the attitude problems and so on. One of the little number of relevant books including the subject is \cite{vlad}. 

For modeling the rigid body dynamics and kinematics in Euler-Poinsot problem the most used canonical coordinates are the Serret-Andoyer (SA) variables. A. Serret discovered them \cite{ser}, his treatment was simplified by Tisserand \cite{tiss} but such approach was diffused by Andoyer (see \cite{and}) who used spherical trigonometry to show how the Serret transformation was only a change of Eulerian coordinates depending upon angular momentum components. Deprit \cite{deprit} established their canonicity. Deprit and Elipe \cite{comp} defined the SA variables not recurring to spherical trigonometry, providing a complete solution of the Euler-Poinsot problem in the phase plane.
The Euler angles reflect the $SO(2)$ symmetry of rotations about the third space axis; but the system is also symmetric with respect to the group $SO(3)$ making possible to reduce the system to only one degree of freedom: this is accomplished by the SA variables.
Their most distinguished feature is then the reduction of the rotational dynamics to a single degree of freedom, because the SA formulation captures the underlying symmetry of the free rigid body problem stemming from conservation of energy and angular momentum. In such a way the single degree of Hamiltonian freedom will yield a phase portrait similar to that of the simple pendulum and with a separatrix confining the librational motion. It can be found e. g. in \cite{tria}.

The SA formalism can be applied to control theory  as in \cite{stab}, where much room is given to the geometry of SA variables and to modeling of control torques. Adlanov and Yudintsev study (see \cite{gyro}) the motion of a spacecraft consisting of a platform whose inertia ellipsoid is triaxial and of a rotor with a small asymmetry with respect to the axis of rotation.
Gurfil (see \cite{basic})  uses the SA modeling of rigid body dynamics in order to obtain some nonlinear stabilizing controllers for an arbitrarily shaped rigid body. Finally, in \cite{kep} Arribas and Elipe study the attitude dynamics of a rigid body on a Keplerian orbit.

\subsection{The herpolhode polar equation by elliptic integrals}

It is well known that a body's torqueless motion around a fixed point O can be characterized with respect to a space-fixed Cartesian coordinate system abstracting from the time law thanks to a geometrical description due to Poinsot. The momental inertia ellipsoid $\mathcal {E}$ referred to O has equation
$Ax^2+By^2+Cy^2=1$. The half-line of the mobile vector $\boldsymbol{\omega}$ intersects such a surface at any time $t$ in a point $Q(t),$ the {\it pole}, which describes on $\mathcal {E}$ an algebraic quartic curve, called the {\it polhode}.

It can be proved that the plane $\mathcal{P}$ tangent to the inertia ellipsoid at $Q(t)$, is fixed in the space. In such a way while the body is revolving about O, simultaneously $\mathcal{E}$ is revolving, maintaining contact at the pole $Q(t)$ with the invariant plane $\mathcal{P}$ which is orthogonal to the fixed direction of the moment of momentum $\bf{K}_O$. As the body moves, the mobile vector $\boldsymbol{\omega}$ traces out the curve {\it herpolhode} in an annulus where the inner (outer)  boundary circle of  center A  corresponds to maximum (minimum) magnitude of $\boldsymbol{\omega}$. The word {\it herpolhode} is formed by three Greek words: $\acute{\epsilon}\rho\pi\epsilon\iota\nu=$ to creep; $\pi \acute{ o}\lambda { o}\varsigma=$ pivot and $ {\rm o}\delta \acute{ o} \varsigma=$ path.
Then it means: \lq\lq pivot serpentine curve''. It is possible to prove that the polhode:
\begin{itemize}
\item  is a closed curve on $\mathcal {E}$, 
\item it is algebraic and not planar, 
\item it may reduce to a point. 
\end{itemize}
In contrast, the herpolhode:
\begin{itemize}
\item could be an open curve, 
\item is transcendental and planar, 
\item is contained within an annulus centered in $Q$ whose bounding circles have radii $\rho_1$ and $\rho_2$,
\item is always concave toward $Q,$ and, when open,  everywhere dense in it.
\end{itemize}
It cannot therefore have points of inflections (see Routh, {\it Dynamics} \cite{rout} p. 472). Furthermore Greenhill (\cite{gree} pp. 227--237) reports that the original herpolhodes drawn by Poinsot (1852) were represented with points of inflection as curves undulating between two concentric circles. But it was pointed out by Hess (1880) and De Sparre (1884) that such points cannot exist. 

All the analytical treatments on the herpolhode are quite similar in their development. We refer to  Appell  \cite{appe2}, \cite{appe3}, citing his conclusion as our starting point. He writes the first integrals of energy and magnitude of moment of momentum as
\begin{equation}\label{DiDi}
\begin{cases}
Ap^2+Bq^2+Cr^2=Dm^2\\
\\
A^2p^2+B^2q^2+C^2r^2=D^2m^2
\end{cases}
\end{equation}
where  $E_0=Dm^2$ and $\|{\bf K}_O\|=D^2m^2$with $m$ homogenizing constant, so that  $D$ has the same nature of $A,\, B,\, C$. After having introduced the further constant
\[
\mathcal{G}=\dfrac{(A-D)(B-D)(C-D)}{ABCD}
\]
he writes, see \cite{appe2} p. 188, formula (39), the differential equation of the herpolhode's  polar curve $\chi=\chi(\rho)$ as:
\[
\rm d\chi=  \dfrac{(\rho^2+\mathcal{G})}{\rho\sqrt{D}\sqrt{-(\rho^2-a)(\rho^2-b)(\rho^2-c)}}\,  \rm d\rho
\]
being
\begin{equation}\label{a+b+c}
\begin{cases}
a=-\dfrac{(B-D)(C-D)}{BCD}\\
\\
b=-\dfrac{(C-D)(A-D)}{CAD}\\
\\
c=-\dfrac{(A-D)(B-D)}{ABD}
\end{cases}
\end{equation}
The vector radius $\rho$ will oscillate between a minimum value  $\rho_{min}=\sqrt{a}$ and  a maximum value $\rho_{max}=\sqrt{b}$. Having assumed $A>B>C$ and $B<D<C$ we have $a>0,\,b>0,\,c<0.$
At p.38, the book \cite{teodorescu} provides a clear representation of the 3D-situation and at p. 322 a planar plot of the herpolhode drawn in a circular annulus.
Then:
\begin{equation}\label{erpoloide}
\chi-\chi_0=\frac{1}{\sqrt{D }} \int_{\rho_{min}}^{\rho}\dfrac{\rho\rm d \rho}{\sqrt{(\rho^2-a)(b-\rho^2)(\rho^2-c)}}+\frac{\mathcal{G}}{\sqrt{D }} \int_{\rho_{min}}^{\rho}\dfrac{\rm d \rho}{\rho\sqrt{(\rho^2-a)(b-\rho^2)(\rho^2-c)}}
\end{equation}
so that one is required to compute two integrals: the first can be driven to a Legendre incomplete integral of the first kind, see the integral $I_5$ of the Appendix. The second is a hyperelliptic one, which can be reduced  to a combination of elliptic integrals of first and third kind, see the integral $I_6$ of the Appendix.
We can therefore read the herpolhode polar equation (i.e. the polar anomaly $\chi$ as a function of ${\rho}$) as:
\begin{equation}\label{finalerp}
\chi(\rho)=\chi_0+\dfrac{1}{2\sqrt{D(a-c)}}\ F\left(\varphi,\ k^{\ast}\right)+\dfrac{\mathcal{G}}{bc\sqrt{D(a-c)}}\left[b\,F(\varphi,k^{\ast})-(b-c)\Pi\left(\varphi,\dfrac{c}{b}{k^{\ast}}^2,k^{\ast}\right)\right]
\end{equation}
where:
\begin{equation}
\begin{cases}
 \varphi = \arcsin{\sqrt{\dfrac{(a-c)(\rho^2-b)}{(a-b)(\rho^2-c)}}}\\
\\
k{\ast}=\sqrt{\dfrac{a-b}{a-c}}
\end{cases}
\end{equation}
Of course, according to the specific problem one wants to solve, the above functions can degenerate into special or elementary ones. 
Whenever the distance between $O$ and $G$ equates the length of a semiaxis of the $\mathcal{E}$, the integration becomes elementary, leading to a curve named a {\it Poinsot spiral} quite carefully studied in \cite{Gomes} pp. 86--89.
For instance, Rosenberg \cite{Rein} mentions a study by Lainé \cite{Lai} who surveyed on a homogeneous lamina shaped as an isosceles triangle which can turn freely about a fixed point. The relevant herpolhode is found by detecting first of all $p,\, q,\, r$, computing the Euler angles, and finally considering the specific locus of the endpoints of $\boldsymbol{\omega}$ over the invariant plane. He obtains a particular law of anomaly $\chi(\rho)$ which coincides with the already mentioned spiral of Poinsot.

Lawden \cite{law} would like do an ample treatment of the subject, but it concludes the study at pp. 135--139, with  parametric formulae providing $\rho$ as a Jacobian elliptic function of time whether the anomaly is given through a complex variable function of time. Almost the same can be read on the mentioned treatises of Whittaker and Appell.

\subsection{The herpolhode equation solved hypergeometrically}
As explained in the Appendix, the integrals $I_5$ and $I_6$ can also be evaluated hypergeometrically. We have:
\begin{equation}\label{finalhyperp}
\chi(\rho)=\chi_0+\dfrac{1}{\sqrt{D}}\sqrt{\frac{\rho^2-b}{(a-b) (b-c)}}\ \mathcal{A}+\dfrac{\mathcal{G}}{\sqrt{D}}\frac{1}{b}\sqrt{\frac{\rho^2-b}{(a-b) (b-c)}}\ \mathcal{L}
\end{equation}
where:
\begin{equation}
\mathcal{A}=\mathrm{F}_{1}\left( \left. 
\begin{array}{c}
\frac12;\frac12,\frac12 \\[2mm]
\frac32
\end{array}
\right|\frac{\rho^2-b}{a-b},-\frac{\rho^2-b}{b-c}\right)
\end{equation}
is the Appell hypergeometric function, despite he does not use it in his treatments \cite{appe2} and \cite{appe3}. In addition $\mathcal{L}$ denotes a Lauricella hypergeometric function:
\begin{equation}
\mathcal{L}=\mathrm{F}_{D}^{(3)}\left( \left. 
\begin{array}{c}
\frac12;1,\frac12,\frac12 \\[2mm]
\frac32
\end{array}
\right|-\frac{\rho^2-b}{b},\frac{\rho^2-b}{a-b},-\frac{\rho^2-b}{b-c}\right)
\end{equation}

\section{The symmetrical body under viscous torques}

The analysis of the heavy body has been performed with the body moving {\it in vacuo} under the sole effect of its own weight which, joint to the constraint reaction, produces a torque whose components are depending, see \eqref{Eu}, on $\gamma$'s, or, generally, on all nine cosines. In both problems we did not meet the Euler equations directly, but worked on their first integrals finding  Euler angles versus time and afterwards the components of instantaneous rotational speed. The case we are now going to analyze keeps the symmetry of the body but without the weight, adding  dissipation: in facts, air resistance acts as an obvious cause of energy loss for the motion of the top. Since it sets the air into motion which is partly communicated to more distant air layers, kinetic energy continually flows from the moving top into the surrounding medium. 
We take into account the effect of viscous torques applied to each of the body axes: the purpose is to detect first how the rotational speed and Euler angles change during time. Klein and Sommerfeld \cite{klei} assumed a drag equivalent to a torque vector with two components: one along the symmetry axis of the body and the other along the orthogonal projection of $\boldsymbol{\omega}$ on the equatorial plane.
We are instead going to model drag effects with respect to all three principal axes of the body, setting the components of the air resistance torques with respect to them as $-\mu\, p,\, -\mu\, q,\,  -\mu r$. Perhaps it would be indicated to choose the coefficient of $r$ smaller than those of $p$ and $q$ since the air will be entrained less by a rotation about the figure axis than by a rotation of the figure axis about an axis perpendicular to it. However, our starting point cannot claim to correspond precisely to the conditions of reality, so we will be satisfied with the {\it common $\mu$ value} approximation and any other frictional influence will be of course disregarded. The authors  of \cite{klei} do not compute how $p$ and $q$ change during time: they merely confine (p. 587) to recognize that the absolute value of $p + iq$ (that is, the length of the equatorial component of the rotation vector, being $i$ the imaginary unit) decreases according to a simple law as $r(t)$, namely:
\[
\sqrt{p^2 +q^2}= \sqrt{p_0^2 +q_0^2} e^{-\mu t/A}.
\]
In chapter X of \cite{Moulton} entitled \lq\lq The damped gyroscope'', the author puts (p.168, eq.29)  the damping couples as $-\mu\ p,\ -\mu\ q,\ 0.$ After a too long sequence of changes of variables, the system of six ordinary differential equations is transformed (p.174, equation 47), into a system of three differential equations whose solutions are expanded as a power series of $\mu/m$, being $m$ a constant. The approach follows an iterative scheme but the whole machinery lacks of compactness and elegance. The problem data are the body specifications $A,\ C,\ \mu $ and the initial values: $p_0,\ q_0,\ r_0, \gamma_{1,0}, \gamma_{2,0},\gamma_{3,0}$. Let us proceed with the integrations.
\subsection{$\boldsymbol{p,\ q, \ r}$ components}
The Euler equations \eqref{Eu} in absence of weight (Poinsot case), with symmetry $A=B$ but under viscous torques on all the body axes become:
\begin{equation}\label{Eudamp1}
\begin{cases}
A{\dot{p}}-Aqr+Cqr=-\mu p\\
A{\dot{q}}+Apr-Cpr=-\mu q\\
C{\dot{r}}=-\mu r\\
p(0)=p_0,\,q(0)=q_0,\,r(0)=r_0\\
\end{cases}
\end{equation}
To integrate \eqref{Eudamp1} first we observe that the first two equations in \eqref{Eudamp1} can be solved for any given function $r(t)$ and such $r(t)$ is immediately computed from the third equation in \eqref{Eudamp1} as 
\begin{equation}\label{exp}
r(t)=r_0 e^{-\mu t/C}.
\end{equation} 
Using a standard calculus technique from the first two equations in \eqref{Eudamp1} it is possible to derive a second order linear differential equation for the unknown $p(t)$. The computations are elementary but tedious, so we chose to omit the details and report the differential equation for $p(t)$
\begin{equation}\label{4.1p}
\ddot{p}+\frac{\mu  (A+2 C) }{A C}\,\dot{p}+\frac{\left[\mu ^2 (A+C)+C r_0^2 (A-C)^2 e^{-\frac{2 \mu  t}{c}}\right]}{A^2C}\,p=0
\end{equation}
Equation \eqref{4.1p} is a linear differential equation of the form
\begin{equation}\label{4.1pol}
\ddot{p}+a\dot{p}+\left(be^{\lambda t}+c\right) p=0
\end{equation}
whose solution can be found in \cite{polyanin} entry 10 at p. 247, expressed by means of Bessel functions (in general of non integer order) of first and second kind $J_\nu$ and $Y_\nu,$ namely:
\begin{equation}\label{4.1pols}
p=e^{-\frac{a}{2}t}\left[C_1J_\nu\left(\frac{2 \sqrt{b}}{\lambda}\,e^{\frac{\lambda  t}{2}}\right)+C_2Y_\nu\left(\frac{2 \sqrt{b}}{\lambda}\,e^{\frac{\lambda  t}{2}}\right)\right]
\end{equation}
being
\begin{equation}\label{4.1nu}
\nu=\frac{\sqrt{a^2-4 c}}{\lambda }
\end{equation}
Going back to \eqref{4.1p}, we see that the parameter setting of our case gives $\nu=-1/2$ so that the Bessel functions assume the form
\[
J_{-\frac12}(x)=\sqrt{\frac{2}{\pi x}}\,\cos x,\quad Y_{-\frac12}(x)=\sqrt{\frac{2}{\pi x}}\,\sin x
\]
In conclusion, the solution to \eqref{4.1p} is given by:
\begin{equation}\label{4.psol}
p(t)=e^{-\frac{\mu\,t}{A}}\left(C_1\, J(t) +C_2\, Y(t)\right)
\end{equation}
Where
\begin{equation}\label{4.psolb}
\begin{split}
J(t)&=\sqrt{\frac{2A \mu }{\pi C r_0 \sqrt{(A-C)^2}}}\, \cos \left[\frac{C r_0 \sqrt{(A-C)^2} e^{-\frac{\mu 
   t}{C}}}{A \mu }\right]\\
   Y(t)&=\sqrt{\frac{2A \mu }{\pi C r_0 \sqrt{(A-C)^2}}}\, \sin \left[\frac{C r_0 \sqrt{(A-C)^2} e^{-\frac{\mu 
   t}{C}}}{A \mu }\right]
   \end{split}
\end{equation}
Notice that in \eqref{4.psolb} we did not simplify the term $\sqrt{(A-C)^2}$ on purpose, since we did not make assumptions regarding $A$ and $C$. Next, the  $\boldsymbol{\omega}$ last  component $q(t)$  comes without integration, since:
\begin{equation}\label{qdit}
q(t)= \frac{e^{\frac{\mu  t}{C}} \left(A \dot{p}+\mu  p(t)\right)}{r_0 (A-C)}.
\end{equation}
The constants $C_1$ and $C_2$ are detected by the initial conditions appearing in \eqref{Eudamp1}, using all the power of Mathematica$^{\text{\textregistered}}$
\begin{equation}
\begin{split}
C_1&=\sqrt{\frac{\pi }{2}} \sqrt{\frac{C r_0 \sqrt{(A-C)^2})}{A \mu }} \left[p_0 \cos \left(\frac{C r_0
   \sqrt{(A-C)^2}}{A \mu }\right)+q_0 \sin \left(\frac{C r_0 \sqrt{(A-C)^2}}{A \mu }\right)\right]\\
   C_2&=\sqrt{\frac{\pi }{2}} \sqrt{\frac{C r_0 \sqrt{(A-C)^2})}{A \mu }}  \left[p_0\,\sin \left(\frac{C r_0 \sqrt{(A-C)^2})}{A \mu
   }\right)-q_0 \cos \left(\frac{C r_0 \sqrt{(A-C)^2})}{A \mu }\right)\right]
   \end{split}
\end{equation}
thus the problem is fully solved. Therefore we gained a closed form expression for each of the rotational speed components  with the motion induced by the sole initial conditions and hindered by viscous drag.

\subsection{$\gamma$ cosines}
Let us go to the Euler angles detection: the system \eqref{pqr} does not seem tractable, so we will come to such angles through the intermediate evaluation of $\gamma$'s.
We could have solved numerically the Poisson \eqref{Poi} system and next refer to \eqref {pqr} for finding the Euler angles. We wish to give a possible variant of the numerical approach in order to provide a better symmetry to all the computations. For the purpose, we will insert all three $(p,\ q,\ r)$  inside the Poisson \eqref{Poi} system which is re-written below:
 \begin{equation}\label{Po}
\begin{cases}
{\dot{\gamma_{1}}}=r\gamma_{2}-q\gamma_{3}\\
{\dot{\gamma_{2}}}=p\gamma_{3}-r\gamma_{1}\\
{\dot{\gamma_{3}}}=q\gamma_{1}-p\gamma_{2}\\
\gamma_1 (0)=\gamma_{1,0}\\
\gamma_2 (0)=\gamma_{2,0}\\
\gamma_3 (0)=\gamma_{3,0}
\end{cases}
\end{equation}
We start with taking the time derivative of the first of \eqref{Po} which holds $\dot{\gamma_2}$ and $\dot{\gamma_3}$; inserting there their expression provided by the second and third equations of \eqref{Po} we obtain a intermediate equation
\begin{equation}\label{gammunobipunto}
\ddot{\gamma_1}=f_1\gamma_1+f_2\gamma_2+f_3\gamma_3\\
\end{equation}
where 
\begin{equation}\label{effe}
\begin{cases}
f_1=\dot{r}+qp\\
f_2=r^2-q^2\\
f_3=-\dot{q}-rp
\end{cases}
\end{equation}
The \eqref{gammunobipunto} is then again derived with respect of time. Inserting there \eqref{Po} we get:
\begin{equation}\label{gammunotripunto}
\dddot{\gamma_1}=g_1\gamma_1+g_2\gamma_2+g_3\gamma_3\\
\end{equation}
where 
\begin{equation}\label{gi}
\begin{cases}
g_1=\dot{f_1}-rf_2+qf_3\\
g_2=\dot{f_2}+rf_1-pf_3\\
g_3=\dot{f_3}-qf_1+pf_2
\end{cases}
\end{equation}
The system of first of \eqref{Po} and \eqref{gammunobipunto} is solved to $\gamma_2$ and $\gamma_3$, obtaining the relations:
\begin{equation}\label{gamma2e3}
\begin{cases}
\gamma_2=\dfrac{f_3\dot{\gamma_1}+q\ddot{\gamma_1}-qf_1\gamma_1}{rf_3-qf_2}\\
\gamma_3=\dfrac{r\ddot{\gamma_1}-rf_1\gamma_1-f_2\dot{\gamma_1}}{rf_3-qf_2}
\end{cases}
\end{equation}
which, inserted in \eqref {gammunotripunto}, provide:
\begin{equation}\label{resolvent}
\begin{cases}
\dddot{\gamma_1}=\ddot{\gamma_1}\left(\dfrac{qg_2+rg_3}{rf_3-qf_2}\right)+\dot{\gamma_1}\left(\dfrac{g_2f_3-f_2g_3}{rf_3-qf_2}\right)+
{\gamma_1}\left(\dfrac{g1-qf_1g_2-rf_1g_3}{rf_3-qf_2}\right)\\
\gamma_1(0)=\gamma_{1,0}\\
\dot{\gamma_1}(0)=r_0 \gamma_{2,0}-q_0\gamma_{3,0}\\
 \ddot{\gamma_1}(0)=\gamma_{1,0}(r_0^2-q_0^2)+\gamma_{2,0}(\dot{r}_0+q_0p_0)+\gamma_{3,0}(\dot{r}_0-q_0p_0)
\end{cases}
\end{equation}
where the third condition comes by merging \eqref{gammunobipunto} and \eqref{gi}. We then get a third order ordinary differential equation in $\gamma_1$ which can only be solved numerically. In any case the remaining cosine $\gamma_2$ stems through the first of \eqref{gamma2e3}, while  $\gamma_3$ is easier seen from the third equation of \eqref{trio}.
\subsection{Euler angles}
The Euler angles are depending on the cosines, then $\theta$ and $\varphi $ can be found without integration by means of \eqref{gammas} as functions of time. On the contrary, the precession $\psi$ will stem (see \eqref{pqr}) for any value of $t$ from a numerical quadrature:
\begin{equation}\label{psi again}
\psi(t)-\psi_0=\int_{0}^t\frac{r-\dot{\varphi}}{\cos\theta}\,{\rm d}t
\end{equation}
Before closing the subject we will provide our comments on how such viscose problem has been treated in relatively recent works. We start with \cite{ivanova}: the rotation of a body around its axis of symmetry is characterized by a turn tensor. At p. 620 the authoress tackles the rotation of a free axisymmetrical rigid body in resisting medium and the interaction between them is simulated by the moment of the linear viscous friction and the tensor of viscous friction. Formula (5.6) provides a vectorial third order non linear ODE in $\boldsymbol{\omega}$ which is faced by looking for it a solution in the form of a series (formulae 5.7, 5.8). Her method is only hinted and the authoress refers to a (Russian) paper of herself.
Nevertheless what mentioned above is enough to understand that it is quite alien to our work.

Grammel in \cite{grammel} deals with the damped top only about technical applications: in such a way it concerns the gyroscope as a damper, namely a device (named {\it Schiffkreisel}) for improving the ship stability against unwanted sea-induced vibrations. The relevant mathematical modeling has then no connection with ours.

Magnus in his book \cite{magnus} at p.155 provides a very short treatment of the viscous top entitled {\it Der symmetrische Kreisel in einem viskosen Medium}. He assumes each axis of a symmetrical body is loaded by both an active and a resistive moment. Of the relevant system of three non linear ODE he first solves the simplest one and -after manipulating the first two- gets an ODE in the complex unknown  ${\omega}^{\ast}$. The transition to the complex field and the reverse have a clear sense. Unfortunately the final smart formula is illusive for being almost all the integrations difficult to be done, and being impossible at all considering the nested ones. If the loading moments $M_k$ are variable during time, no-one of the written integrals can be solved in closed form.

\section{Conclusions}
 Someone could believe that with the current possibilities of machine computations, the achievement of analytical solutions to rigid body motions should be of academic interest only. It must then be countered that numerical methods enable us to obtain accurate solutions in short time intervals of integration but they {\it cannot capture the long-time behavior of solutions}. The authors of \cite{calvo} refer {about} recent applications such as the attitude evolution of a spinning spacecraft which have opened new requirements on the computer algorithms for onboard computations: so that the main task is not the accuracy but the reliability of the algorithms. Such analytical solutions, keeping in mind that with large angle rotations the theory is non-linear, are then useful in parametric studies, error analyses, onboard computations and stability analysis. Some authors in \cite{long} repeat the same concept as it concerns thrusting/spinning spacecraft problems, but we disagree with their idea that complex variables could significantly contribute {\it today} to the compactness of the final solutions. Not all analytical solutions attain the same level of generalization. E.g. the analytical solutions produced, 1996, by Tsiotras \& Longuski \cite{Tsi} have  heavy restrictions: they do a first order correction to their previous linear-zero-order solution which required a symmetry or near-symmetry assumption and had been accomplished by means of the Fresnel complex integral. The mentioned correction resorts to the exponential integral. Finally, they provide an explicit formula for the bound of the error of the approximation.

 Our analytical solutions, built on special functions of Mathematical Physics, are reasonably general. Let us provide a table with the main outcomes of the paper: those that are new, as far as we are concerned, are:
 \eqref{psilaur}, \eqref{psitorqueless}, \eqref{philaur},  \eqref{resolvent}, \eqref{pqerre}, \eqref{4.psol},\eqref{qdit}, \eqref{finalerp}, \eqref{finalhyperp}.

\begin{center} 
\begin{tabular}{ |l|l|l|l|l| }
  \hline
  \multicolumn{5}{|c|}{\textbf{Problem}} \\
  \hline
 \hline  & \textbf{Lagrange-Poisson} &  \textbf{Euler-Poinsot} & \textbf{Viscous drag} & \textbf{Herpolhode} \\
\hline  Simmetry & $A=B$ & none, $A>B>C$ & $A=B$ & none, $A>B>C$ \\
\hline  Load & weight & none & symmetric drag & none  \\
\hline  Drag & none & none & viscous & none \\
\hline $\theta$ &\eqref{sn} &\eqref{2angoli} &\eqref{gammas} &--\\
\hline  $\psi$ &\eqref{psilaur} &\eqref{psitorqueless}  &\eqref{psi again}  &-- \\
\hline  $\varphi$  &\eqref{philaur}  &\eqref{2angoli}  &\eqref{gammas}  &-- \\
\hline  $\gamma_1$ &\eqref{gammas}  & \eqref{gammas} &\eqref{gammunotripunto}  & --\\
\hline  $\gamma_2$ & \eqref{gammas}&\eqref{gammas} &\eqref{gamma2e3}  & --\\
\hline  $\gamma_3$ &\eqref{gammas} & \eqref{gammas}& \eqref{gamma2e3} &-- \\
\hline  $p$ &\eqref{pq} &\eqref{pqerre} & \eqref{4.psol}& --\\
\hline  $q$ &\eqref{pq}  &\eqref{pqerre}  &\eqref{qdit}  & --\\
\hline  $r$ & \eqref{Quarto} &\eqref{pqerre}  &\eqref{exp}  & --\\
\hline  $\chi$ & --&-- &-- &\eqref{finalerp},\eqref{finalhyperp}\\
  \hline
\end{tabular}
\vspace{2mm}

 Motion about a fixed point: a summary of problems and solutions
\end{center}
 To the above table we add the summary below.

{\it Heavy body}: we obtained an explicit expression of the precession angle through the nutation, being this last known by Jacobian cosinus amplitude of time. The precession and the intrinsic rotation have been found as functions of time, both via the Hypergeometric functions. In such a way the time history of each point of the body becomes computable via the matrix rotation.

{\it Asymmetrical body free of torques}: to the well-known Jacobi formul\ae\, of nutation and intrinsic rotation angles, we add that of precession by means of elliptic integrals of third kind. 

{\it Herpolhode}: we start from its Appell (but really due to Poinsot) formulation getting it in finite terms by means of elliptic integrals of first and third kind or, alternatively,  through Lauricella and Appell hypergeometric functions.

We faced the Euler equations directly when compelled by dissipation: with a {\it symmetrical body under viscous drag}, we find the components $p$ and $q$ of rotational speed as damped goniometric functions of time. The way has then pointed out for computing the $\gamma$ cosines of the Cartesian coordinate system body; then all three Euler angles can be found.

\appendix

\section{Appendix: Hypergeometric identities}\label{identities}

We recall hereinafter an outline of the Lauricella
hypergeometric functions. The first hypergeometric historical series appeared in the Wallis's \textit{Arithmetica
infinitorum} (1656): 
\[
_2\mathrm{F}_{1}\left(\left. 
\begin{array}{c}
a,b\\[2mm]
c
\end{array}
\right|x\right)=1+\frac{a\, b}{1\cdot c}x+\frac{a\, (a+1)\, b%
\, (b+1)}{1\cdot 2\cdot c\cdot (c+1)}x^{2}+\cdots ,
\]
for $|x|<1$ and real parameters $a,\,b,\,c.$ The product of $n$ ascending factors: 
\[
(\lambda )_{n}=\lambda \left( \lambda +1\right) \cdots \left( \lambda
+n-1\right)=\frac{\Gamma(\lambda+n)}{\Gamma(\lambda)},
\]
called \textit{Pochhammer symbol} (or \textit{truncated factorial}
) allows to represent $_{2}F_{1}$ as: 
\[
_2\mathrm{F}_{1}\left(\left. 
\begin{array}{c}
a,b\\[2mm]
c
\end{array}
\right|x\right)=\sum_{n=0}^{\infty }\frac{\left( a\right) _{n}\left(
b\right) _{n}}{\left( c\right) _{n}}\frac{x^{n}}{n!}.
\]
A meaningful contribution on various $_{2}F_{1}$ topics
is ascribed to Euler in three papers \cite{Eul1738, Eul1769, Eul1792}
but he does not seem \cite{dutka} to have known the integral representation: 
\[
_2\mathrm{F}_{1}\left(\left. 
\begin{array}{c}
a,b\\[2mm]
c
\end{array}
\right|x\right)=\frac{\Gamma (c)}{\Gamma (a)\Gamma (c-a)}\,\int_{0}^{1}%
\frac{u^{a-1}(1-u)^{c-a-1}}{(1-xu)^{b}}\,\mathrm{d}u,
\]
really due to Legendre, \cite{legint}, sect. 2. The above integral relationship is true if $c>a>0$ and for 
$\left| x\right| <1,$ even if this limitation can be discarded thanks to the analytic
continuation.\ 

Many functions have been introduced in 19$^{\mathrm{th}}$ century for
generalizing the hypergeometric functions to multiple variables. We recall the Appell ${\rm F}_1$ two--variable hypergeometric series, defined as: 
\[
\mathrm{F}_{1}\left( \left. 
\begin{array}{c}
a;b_{1},b_{2} \\[2mm]
c
\end{array}
\right| x_{1},x_{2}\right) =\sum_{m_{1}=0}^{\infty }\sum_{m_{2}=0}^{\infty }%
\frac{(a)_{m_{1}+m_{2}}(b_{1})_{m_{1}}(b_{2})_{m_{2}}}{(c)_{m_{1}+m_{2}}}%
\frac{x_{1}^{m_{1}}}{m_{1}!}\frac{x_{2}^{m_{2}}}{m_{2}!},\quad |x_{1}|<1,\,|x_{2}|<1.
\]
The analytic continuation of  Appell's function on
$\mathbb{C}\setminus [1,\infty)\times\mathbb{C}\setminus [1,\infty)$ comes from its integral representation theorem: if $\mathrm{Re}\,(a)>0$,\,$\mathrm{Re}\,(c-a)>0$, then:

\begin{equation}
\mathrm{F}_{1}\left( \left. 
\begin{array}{c}
a;b_{1},b_{2} \\[2mm]
c
\end{array}
\right| x_{1},x_{2}\right) =\frac{\Gamma (c)}{\Gamma (a)\Gamma (c-a)}%
\int_{0}^{1}\frac{u^{a-1}\left( 1-u\right) ^{c-a-1}}{\left(
1-x_{1}\,u\right) ^{b_{1}}\left( 1-x_{2}\,u\right) ^{b_{2}}}\,\mathrm{d}u.
\label{F1}
\end{equation}
The functions introduced and investigated by G. Lauricella (1893) and S. Saran (1954), are of our prevailing interest; and among them the hypergeometric
function $F_{D}^{(n)}$ of $n\in \mathbb{N}^{+}$ variables (and $n+2$
parameters), see \cite{s} and \cite{l}, defined as: 
\[
\mathrm{F}_{D}^{(n)}\left(\left. 
\begin{array}{c}
a,b_{1},\ldots ,b_{n}\\[2mm]
c
\end{array}
\right|x_{1},\ldots ,x_{n}\right):=
\sum_{m_{1},\ldots ,m_{n}\in \mathbb{N}}\frac{(a)_{m_{1}+\cdots
+m_{n}}(b_{1})_{m_{1}}\cdots (b_{n})_{m_{n}}}{(c)_{m_{1}+\cdots
+m_{n}}m_{1}!\cdots m_{n}!}\,x_{1}^{m_{1}}\cdots x_{n}^{m_{m}} 
\]
with the hypergeometric series usual convergence requirements $%
|x_{1}|<1,\ldots ,|x_{n}|<1$. If $\mathrm{Re}\,c>\mathrm{Re}\,a>0$ , the
relevant Integral Representation Theorem provides: 
\begin{equation}\label{iirtt}
\mathrm{F}_{D}^{(n)}\left(\left. 
\begin{array}{c}
a,b_{1},\ldots ,b_{n}\\[2mm]
c
\end{array}
\right|x_{1},\ldots ,x_{n}\right)=\frac{%
\Gamma (c)}{\Gamma (a)\,\Gamma (c-a)}\,\int_{0}^{1}\,\frac{%
u^{a-1}(1-u)^{c-a-1}}{(1-x_{1}u)^{b_{1}}\cdots (1-x_{n}u)^{b_{n}}}\,\mathrm{d%
}u 
\end{equation}
allowing the analytic continuation to $\mathbb{C}^{n}$ deprived of the
cartesian $n$-dimensional product of the interval $]1,\infty [$ with itself. 

Notice that the integral formula \eqref{iirtt} allows an easy computer algebra implementation, being \eqref{iirtt} a parametric integral, which is suitable for each possible numerical computation.

\section{Some integrals evaluation}
The reader can find hereinafter some details about the integrations used throughout the text.
\subsection{Integral $I_1$}

Let $a>b>y>c$ and $\alpha$ a real number. Consider
\[
I_1(a,b,c;\alpha;y)=\int_{c}^y\frac{1-\alpha u}{1-u^2}\frac{{\rm d}u}{\sqrt{(a-u)(b-u)(u-c)}}
\]
To evaluate the integral we normalize the interval of integration, putting:
\[
s=\frac{u-c}{y-c}
\]
this leads to this expression for the given integral, where we remark that in the application of our interest we can limit to $0<c<1$:
\[
I_1(a,b,c;\alpha;y)=\frac{\sqrt{y-c}}{(1-c^2)
   \sqrt{(a-c) (b-c)}}\int_0^1\frac{1-c\alpha-\alpha (y-c) s}{\left(1-\frac{y-c}{1-c}s\right) \left(1+\frac{y-c}{1+c}s\right)\sqrt{s\left(1+\frac{
   y-c}{c-a}s\right)\left(1+\frac{y-c}{c-b}s\right)}}{\rm d}s
\]
which is therefore expressible, via the integral representation theorem, in terms of two Lauricella functions of four variables, that is:
\[
I_1(a,b,c;\alpha;y)=\frac{\sqrt{y-c}}{(1-c^2)
   \sqrt{(a-c) (b-c)}}\left(2(1-c\alpha)X-\frac23\alpha(y-c)Y\right)
\]
where
\[
\begin{split}
X&=\mathrm{F}_{D}^{(4)}\left( \left. 
\begin{array}{c}
\frac12;1,1,\frac12,\frac12 \\[2mm]
\frac32
\end{array}
\right| \frac{y-c}{1-c},-\frac{y-c}{1+c},-\frac{y-c}{c-a},-\frac{y-c}{c-b}\right)\\
Y&=\mathrm{F}_{D}^{(4)}\left( \left. 
\begin{array}{c}
\frac32;1,1,\frac12,\frac12 \\[2mm]
\frac52
\end{array}
\right| \frac{y-c}{1-c},-\frac{y-c}{1+c},-\frac{y-c}{c-a},-\frac{y-c}{c-b}\right)
\end{split}
\]
Moreover integral $I_1$ could also have been evaluated in terms of elliptic integrals of third kind, in fact following \cite{byrd} entry 233.20 p. 74 we are lead to
\begin{equation}\label{by233:20}
I_1(a,b,c;\alpha;y)=\frac{1}{\sqrt{a-c}}\left(\frac{1-\alpha }{1-c} X_1+\frac{1+\alpha}{1+c}X_2\right)
\end{equation}
where
\begin{equation}\label{t1-2}
\begin{split}
X_1&=\Pi \left(\text{am}\left(F\left(\arcsin\left(\sqrt{\frac{y-c}{b-c}}\right),\sqrt{\frac{b-c}{a-c}}\right),\sqrt{\frac{b-c}{a-c}}\right)
   ,\frac{b-c}{1-c},\sqrt{\frac{b-c}{a-c}}\right)\\
X_2&=\Pi
   \left(\text{am}\left(F\left(\arcsin
   \left(\sqrt{\frac{y-c}{b-c}}\right),\sqrt{\frac{b-c}{a-c}}\right),\sqrt{\frac{b-c}{a-c}}\right),-\frac{b-c}{c+1},
   \sqrt{\frac{b-c}{a-c}}\right)
   \end{split}
\end{equation}
\subsection{Integral $I_2$}
With the same assumptions of integral $I_1$, we consider
\[
I_2(a,b,c;\alpha;y)=\int_{c}^y\frac{1-\alpha u}{1-u^2}\frac{u}{\sqrt{(a-u)(b-u)(u-c)}}\,{\rm d}u
\]
Using the same change of variable used for $I_1$, we can also express $I_2$ in terms of Lauricella functions, namely:
\[
I_2(a,b,c;\alpha;y)=\frac{\sqrt{y-c}}{(1-c^2)
   \sqrt{(a-c) (b-c)}}\left(2 c (1-c\alpha)X+\frac{2}{3} (1-2 c \alpha) (y-c)Y-\frac{2}{5}\alpha  \left((c-y)^2\right)Z\right)
\]
where
\[
\begin{split}
X&=\mathrm{F}_{D}^{(4)}\left( \left. 
\begin{array}{c}
\frac12;1,1,\frac12,\frac12 \\[2mm]
\frac32
\end{array}
\right| \frac{y-c}{1-c},-\frac{y-c}{1+c},-\frac{y-c}{c-a},-\frac{y-c}{c-b}\right)\\
Y&=\mathrm{F}_{D}^{(4)}\left( \left. 
\begin{array}{c}
\frac32;1,1,\frac12,\frac12 \\[2mm]
\frac52
\end{array}
\right| \frac{y-c}{1-c},-\frac{y-c}{1+c},-\frac{y-c}{c-a},-\frac{y-c}{c-b}\right)\\
Z&=\mathrm{F}_{D}^{(4)}\left( \left. 
\begin{array}{c}
\frac52;1,1,\frac12,\frac12 \\[2mm]
\frac72
\end{array}
\right| \frac{y-c}{1-c},-\frac{y-c}{1+c},-\frac{y-c}{c-a},-\frac{y-c}{c-b}\right)
\end{split}
\]
Again following \cite{byrd} entry 233.20 p. 74 we can express $I_2$ in terms of elliptic integral of third and first kind as
\begin{equation}\label{by233:20b}
I_2(a,b,c;\alpha;y)=\frac{1}{\sqrt{a-c}}\left(\frac{1-\alpha  }{1-c} X_1-\frac{1+\alpha }{1+c}X_2+2\alpha X_3\right)
\end{equation}
where $X_1$ and $X_2$ are the same as \eqref{t1-2} while
\begin{equation}\label{}
X_3=F\left(\arcsin
   \left(\sqrt{\frac{y-c}{b-c}}\right),\sqrt{\frac{b-c}{a-c}}\right)
\end{equation}

\subsection{Integral $I_3$}
For $c>1$ we consider
\[
I_3(y,c)=\int_0^y\frac{{\rm sn}^2u}{1-c^2{\rm sn}^2u}\,{\rm d}u
\]
The integral appears as entry 337.01 p. 201 of \cite{byrd}, namely
\[
I_3(y,c)=\frac{1}{c^2}\left[\Pi({\rm am}(y),c^2,k)-y\right]
\]
\subsection{Integral $I_4$}
For $c<1$ we consider:
\[
I_4(y,c)=\int_0^y\frac{{\rm cn}^2u}{1-c^2{\rm sn}^2u}\,{\rm d}u
\]
The integral appears as entry 338.01 p. 202 of \cite{byrd}, namely
\[
I_4(y,c)=\frac{1}{c^2}\left[(c^2-1)\Pi({\rm am}(y),c^2,k)+y\right]
\]
\subsection{Integral $I_5$}
If $c<a<y<b$ consider the integral
\[
I_5(a,b,c;y)=\int_{\sqrt{b}}^y\frac{u}{\sqrt{(a-u^2)(u^2-b)(u^2-c)}}\,{\rm d}u
\]
A natural step is a change of variable $u^2=x$ which gives the elliptic integral of first kind
\begin{equation}\label{intmed1}
I_5(a,b,c;y)=\frac12\int_{b}^{y^2}\frac{{\rm d}x}{\sqrt{(a-x)(x-b)(x-c)}}
\end{equation}
which is computed in \cite{grad} entry 2.131.5 p. 250 and \cite{byrd} entry 235.00 leading to:
\[
I_5(a,b,c;y)=\frac{1}{\sqrt{a-c}}F(\varphi,k)
\]
being
\[
\varphi=\arcsin\sqrt{\frac{(a-c)(y^2-b)}{(a-b)(y^2-c)}},\quad k^2=\frac{a-b}{a-c}
\]
\subsubsection{Remark}
This integral can, of course, be evaluated using the hypergeometric approach. In fact, if in \eqref{intmed1} we introduce the change of variable $s=\frac{x-b}{y^2-b}$ we obtain:
\[
I_5(a,b,c;y)=\frac12\sqrt{\frac{y^2-b}{(a-b) (b-c)}}\int_0^1\frac{s^{-\frac12}}{  \sqrt{\left(1-\frac{
   y^2-b}{a-b}s\right) \left(1+\frac{y^2-b}{b-c}s\right)}}\,{\rm d}s
\]
Therefore we can use the relevant integral representation, which allows us to express $I_5$ in terms of the Appell two variable hypergeometric function:
\[
I_5(a,b,c;y)=\sqrt{\frac{y^2-b}{(a-b) (b-c)}}\;\mathrm{F}_{1}\left( \left. 
\begin{array}{c}
\frac12;\frac12,\frac12 \\[2mm]
\frac32
\end{array}
\right|\frac{y^2-b}{a-b},-\frac{y^2-b}{b-c}\right)
\]
\subsection{Integral $I_6$}
With the same parameters of $I_5$ we consider
\[
I_6(a,b,c;y)=\int_{\sqrt{b}}^y\frac{{\rm d}u}{u\sqrt{(a-u^2)(u^2-b)(u^2-c)}}
\]
Again, the change of variable $u^2=x$ leads to an elliptic integral, namely
\begin{equation}\label{intmed2}
I_6(a,b,c;y)=\frac12\int_{b}^{y^2}\frac{{\rm d}x}{x\sqrt{(a-x)(x-b)(x-c)}}
\end{equation}
This last integral is tabulated in \cite{grad} entry 3.137.5 p. 259 and in \cite{byrd} entry 235.17, leading to:
\[
I_6(a,b,c;y)=\frac{1}{bc\sqrt{a-c}}\left[b\,F(\varphi,k)-(b-c)\Pi\left(\varphi,\frac{c}{b}k^2,k\right)\right]
\]
where $\varphi$ and $k$ are the same introduced for the fifth integral.
\subsubsection{Remark}
This integral also can be evaluated using the hypergeometric approach, so using for \eqref{intmed2} the same change of variable, used for the fifth integral, we obtain the expression:
\[
I_6(a,b,c;y)=\frac{1}{2b}\sqrt{\frac{y^2-b}{(a-b) (b-c)}}\int_0^1\frac{s^{-\frac12}}{\left(1+\frac{y^2-b}{b}s \right)  \sqrt{\left(1-\frac{
   y^2-b}{a-b}s\right) \left(1+\frac{y^2-b}{b-c}s\right)}}\,{\rm d}s
\]
which drives, using the integral representation, to express $I_6$ in terms of a Lauricella function of three variables:
\[
I_6(a,b,c;y)=\frac{1}{b}\sqrt{\frac{y^2-b}{(a-b) (b-c)}}\;\mathrm{F}_{D}^{(3)}\left( \left. 
\begin{array}{c}
\frac12;1,\frac12,\frac12 \\[2mm]
\frac32
\end{array}
\right|-\frac{y^2-b}{b},\frac{y^2-b}{a-b},-\frac{y^2-b}{b-c}\right)
\]
\subsection{Integral $I_7$}

We go back to the parameters of the first and second integral $a>b>y>c$. Here the integral is quite a simple one, but here we are concerned about its inversion. The integral, given in \cite{grad} entry 3.131.3 p. 230 and \cite{byrd} entry 233.00 is:
\[
I_7(a,b,c;y)=\int_{c}^y\frac{{\rm d}u}{\sqrt{(a-u)(b-u)(u-c)}}=\frac{2}{\sqrt{a-c}}\,F\left(\arcsin\sqrt{\frac{y-c}{b-c}},\sqrt{\frac{b-c}{a-c}}\right)
\]
The integral inversion, that is the solution with respect to $y$ of the equation $I_7(a,b,c;y)=L$ is obtained by recalling the Jacobi amplitude ${\rm am}(u,k),$ which is the inverse of the elliptic integral of first kind $F({\rm am}(u,k),u)=u$ and the Jacobi sinus amplitude ${\rm sn}(u,k)=\sin{\rm am}(u,k).$ The inversion formula is
\[
I_7(a,b,c;y)=L\iff y=c+(b-c)\, \text{sn}^2\left(\frac{L\,\sqrt{a-c}}{2} ,\sqrt{\frac{b-c}{a-c}}\right)
\]
Observe that the equation has a solution if $L$ is such that
\[
L\leq\frac{2}{\sqrt{a-c}}\,{\bf K}\left(\sqrt{\frac{b-c}{a-c}}\right)
\]

After having cited some books about Lauricella functions we deem to quote papers where their application is best shown as \cite{msr} and \cite{Kr}.

\subsubsection*{Acknowledgments}

The second author is partially supported by an Italian RFO research grant.

\end{document}